\documentclass[12pt,a4paper]{article}
\usepackage{amsmath,amssymb}
\usepackage{diagrams}
\usepackage[thmmarks]{ntheorem}
\usepackage{tikz}

\newcommand{\oAss}{\mathcal{A}\mathit{ss}}
\newcommand{\oA}{\mathcal{A}}
\newcommand{\oB}{\mathcal{B}}
\newcommand{\oC}{\mathcal{C}}	

\newcommand{\oM}{\mathcal{M}}
\newcommand{\oP}{\mathcal{P}}
\newcommand{\oR}{\mathcal{R}}
\newcommand{\oEnd}{{\mathcal{E}\mathit{nd}}}

\newcommand{\alg}{\alpha}	

\newcommand{\ddA}{\dd_\oA}	

\newcommand{\ddR}{{\dd_\oR}}	
\newcommand{\rR}{{\rho_\oR}}	

\newcommand{\oDA}{\mathcal{DA}}	
\newcommand{\ddDA}{\dd_\oDA}	
\newcommand{\oDR}{\mathcal{DR}}	
\newcommand{\ddDR}{{\dd_\oDR}}	
\newcommand{\rDR}{{\rho_\oDR}}	

\newcommand{\olDR}{{\ol{{\oDR}}}}	
\newcommand{\ddolDR}{\dd_{\ol{\oDR}}}	
\newcommand{\interm}{\iota}	
\newcommand{\rolDR}{{\rho_{\olDR}}}	

\newcommand{\oMDA}{\mathcal{MDA}}	
\newcommand{\ddMDA}{\dd_\oMDA}	
\newcommand{\olMDR}{\mathcal{M\ol{DR}}}	
\newcommand{\ddolMDR}{{\dd_{\olMDR}}}	
\newcommand{\rolMDR}{{\rho_{\olMDR}}}	
\newcommand{\oMR}{\mathcal{MR}}	

\newcommand{\set}[2]{\left\{ #1 \ | \ #2 \right\} }	
\newcommand{\coll}[1]{\ensuremath{#1}-collection}	
\newcommand{\colop}[1]{\ensuremath{#1}-operad}	
\newcommand{\ar}[1]{{\mathrm{ar}(#1)}}
\newcommand{\gr}{{\mathrm{gr}}}
\newcommand{\dd}{\partial}	
\newcommand{\codd}{\delta}	
\newcommand{\coddH}{\codd_{\mathrm{Hoch}}}	
\newcommand{\dg}[1]{{\left| #1 \right|}}	
\newcommand{\Fr}[2]{ \mathbb{F}^{ #1 }\left( #2 \right) } 
\newcommand{\FrM}[2]{ #1 \left\langle #2 \right\rangle }	
\newcommand{\id}{\mathsf{1}}	
\newcommand{\into}{\hookrightarrow}	
\newcommand{\onto}{\twoheadrightarrow}	
\newcommand{\inp}[1]{I(#1)}	
\newcommand{\out}[1]{O(#1)}	
\newcommand{\N}{\mathbb{N}}
\newcommand{\Z}{\mathbb{Z}}
\newcommand{\cC}{\mathsf{C}}	
\newcommand{\Phitot}{\Phi_{\Sigma}}	
\newcommand{\BC}{\Omega\mathrm{B}}	

\newcommand{\upar}{\raisebox{1pt}[0pt][0pt]{\ensuremath{\uparrow}}}	
\newcommand{\downar}{\raisebox{1pt}[0pt][0pt]{\ensuremath{\downarrow}}}	
\newcommand{\kspan}[1]{k\!\left\langle #1 \right\rangle}	
\newcommand{\ul}[1]{\underline{#1}}	
\newcommand{\ol}[1]{\overline{#1}}	
\newcommand{\op}{\oplus}	
\newcommand{\ot}{\otimes}	
\newcommand{\otexp}[2]{{#1^{\ot #2}}}	
\newcommand{\oo}{\circ}	

\newcommand{\proj}{\mathop{\mathrm{pr}}\nolimits}	
\newcommand{\Ker}{\mathop{\mathrm{Ker}}\nolimits} 
\renewcommand{\Im}{\mathop{\mathrm{Im}}\nolimits}	
\newcommand{\fpr}{\mathop{\mathrm{*}}\nolimits}	
\newcommand{\bigfpr}{\raisebox{-6pt}[0pt][0pt]{\mbox{\Large *}}}	
\newcommand{\pre}[1]{\mathop{\mathrm{pre}(#1)}\nolimits}	
\newcommand{\dgmod}[1]{\mathrm{dg}-#1-\mathrm{mod}}	
\newcommand{\Der}{\mathop{\mathrm{Der}}\nolimits} 
\newcommand{\colim}{\mathop{\mathrm{colim}}}	
\newcommand{\Ext}{\mathop{\mathrm{Ext}}\nolimits}
\newcommand{\Hom}{\mathop{\mathrm{Hom}}\nolimits}

{
\theoremstyle{change}
\theorembodyfont{\normalfont}
\theoremseparator{.}

\newtheorem{definition}{Definition}[section]
}

{
\theoremstyle{change}
\theorembodyfont{\normalfont}
\theoremseparator{.}

\newtheorem{theorem}[definition]{Theorem}
}

{
\theoremstyle{change}
\theorembodyfont{\normalfont}
\theoremseparator{.}

\newtheorem{corollary}[definition]{Corollary}
}

{
\theoremstyle{change}
\theorembodyfont{\normalfont}
\theoremseparator{.}

\newtheorem{lemma}[definition]{Lemma}
}

{
\theoremstyle{change}
\theorembodyfont{\normalfont}
\theoremseparator{.}

\newtheorem{sublemma}[definition]{Sublemma}
}

{
\theoremstyle{change}
\theorembodyfont{\normalfont}
\theoremseparator{.}

\newtheorem{example}[definition]{Example}
}

{
\theoremstyle{change}
\theorembodyfont{\normalfont}
\theoremseparator{.}

}

{
\theoremstyle{change}
\theorembodyfont{\normalfont}
\theoremseparator{.}

\newtheorem{convention}[definition]{Convention}
}

{
\theoremstyle{nonumberplain}
\theorembodyfont{\normalfont}
\theoremseparator{.}

}

{
\theoremstyle{nonumberplain}
\theoremheaderfont{\it}
\theoremseparator{.}
\theorembodyfont{\normalfont}
\theoremsymbol{\ensuremath{\square}}

\newtheorem{proof}{Proof}
}

{
\theoremstyle{change}
\theorembodyfont{\normalfont}

}

\title{Gerstenhaber-Schack diagram cohomology from operadic point of view}
\author{
Martin Doubek\thanks{The author was supported by GA\v{C}R 201/09/H012.},  \\
Charles University, Prague\\
\texttt{martindoubek@seznam.cz}
}
	
\date{\today}

\begin{document}

\qedsymbol{\ensuremath{\square}}

\maketitle

\abstract{
We show that the operadic cohomology for any type of algebras over a non-symmetric operad $\oA$ can be computed as $\Ext$ in the category of operadic $\oA$-modules.
We use this principle to prove that the Gerstenhaber-Schack diagram cohomology is operadic cohomology.
}

\section{Introduction}

The Operadic Cohomology (OC) gives a systematic way of constructing cohomology theories for algebras $A$ over an operad $\oA$.
It recovers the classical cases: Hochschild, Chevalley-Eilenberg, Harrison etc.
It also applies to algebras over coloured operads (e.g. morphism of algebras) and over PROPs (e.g. bialgebras).
The OC first appeared in papers \cite{Cotangent}, \cite{MMBialg} by M. Markl.

Abstractly, the OC is isomorphic to the triple cohomology, at least for algebras over Koszul operads \cite{FMDL}.
It is also isomorphic to the Andr\'e-Quillen Cohomology (AQC).
In fact, the definition of OC is analogous to that of AQC:
It computes the derived functor of the functor $\Der$ of derivations like AQC, but does so in the category of operads.
While AQC offers a wider freedom for the choice of a resolution of the given algebra $A$, OC uses a particular universal resolution for all $\oA$-algebras (this resolution is implicit, technically OC resolves the \emph{operad} $\oA$).
Thus there is, for example, a universal construction of an $L_\infty$ structure on the complex computing OC \cite{IB}, whose generalized Maurer-Cartan equation describes formal deformations of $A$.

The success of OC is due to the Koszul duality theory \cite{GK}, which allows us to construct resolutions of Koszul operads explicitly.
Koszul theory has received a lot of attention recently \cite{Vallette} and now goes beyond operads.
However, it still has its limitations:

On one hand, it is bound to quadratic relations in a presentation of the operad $\oA$.
The problem with higher relations can be remedied by using a different presentation, but it comes at the cost of increasing the size of the resolution (e.g. \cite{ValletteBV}). 
This is not a major problem in applications, but the minimal resolutions have some nice properties - namely they are unique up to an isomorphism thus providing a cohomology theory unique already at the chain level.
So the construction of the minimal resolutions is still of interest.

On the other hand, there are quadratic operads which are not Koszul and for those very little is known \cite{MR}.

In this paper, we show that OC is isomorphic to $\Ext$ in the category of operadic $\oA$-modules.
Thus instead of resolving the operad $\oA$, it suffices to find a projective resolution of a specific $\oA$-module $\oMDA$ associated to $\oA$. 
The ideas used here were already sketched in the paper \cite{IB} by M. Markl.

The resolution of operadic modules are probably much easier to construct explicitly than resolutions of operads, though this has to be explored yet.
This simplification allows us to make a small step beyond Koszul theory:

An interesting example of a non-Koszul operad is the coloured operad describing a diagram of a fixed shape consisting of algebras over a fixed operad and morphisms of those algebras.
The case of a single morphism between two algebras over a Koszul operad is long well understood.
For a morphism between algebras over a general operad as well as for diagrams of a few simple shapes, some partial results were obtained in \cite{HDA}.
These are however not explicit enough to write down the OC.

On the other hand, a satisfactory cohomology for diagrams was invented by Gerstenhaber and Schack \cite{GS} in an ad-hoc manner.
In \cite{FMY}, the authors proved that the Gerstenhaber-Schack cohomology of a single morphism of associative or Lie algebras is operadic cohomology.
We use our theory to extend this result to arbitrary diagrams.

The method used can probably be applied in a more general context to show that a given cohomology theory is isomorphic to OC.
The original example is \cite{IB} (and similar approach also appears in \cite{MV}), where the author proves that Gerstenhaber-Schack \emph{bialgebra} cohomology is the operadic cohomology.
Also the method might give an insight into the structure of operadic resolutions themselves, the problem we won't mention in this paper.

On the way, we obtain a modification of the usual OC which includes the quotient by infinitesimal automorphisms (Section \ref{AugmentedCotangentComplex}).

Also an explicit description of a free resolution of the operad $\oA$ \emph{with adjoined derivation} is given if a free resolution of $\oA$ is explicitly given.
This appeared already in \cite{IB} and produces several new examples of minimal resolutions and as such might be of an independent interest.

We assume the reader is familiar with the language of operads (e.g. \cite{MSS},\cite{LV}).

Finally, I would like to thank Martin Markl for many useful discussions.

\begin{convention}
As our main object of interest is a diagram of associative algebras, we will get by with \emph{non-symmetric} operads, that is operads with no action of the permutation groups.
The results can probably be generalized in a straightforward way to symmetric operads.
\end{convention}

\pagebreak

\tableofcontents
\bigskip

In \textbf{Section \ref{SectionBasics}}, we briefly recall basic notions of the operad theory with focus on \emph{coloured} operads (see also \cite{HAAHA},\cite{HDA}).
We pay special attention to operadic modules.
We introduce the notion of tree composition which is just a convenient way to write down complicated operadic compositions.
In Section \ref{SectionFreeProduct}, we discuss free product of operads and obtain a form of the K\"unneth formula computing homology of the free product.

In \textbf{Section \ref{SectionOC}}, we develop the theory sketched by M. Markl in Appendix B of \cite{IB}.
We give full details for coloured operads.
We begin by recalling the operadic cohomology.
In Section \ref{AlgebrasWithDerivation}, we construct an explicit resolution of the operad $\oDA$ describing algebras over $\oA$ with adjoined derivation assuming we know an explicit resolution of the operad $\oA$.
In Section \ref{AugmentedCotangentComplex}, we clarify the significance of operadic derivations on the resolution of $\oDA$ with values in $\oEnd_A$.
This leads to an augmentation of the cotangent complex which has nice interpretation in terms of formal deformation theory.
In Sections \ref{SectionIntermediate} and \ref{FromOperadsToModules}, we realize that all the information needed to construct augmented cohomology is contained in a certain operadic \emph{module}.
This module is intrinsically characterized by being a resolution (in the category of operadic modules) of $\oMDA$, a certain module constructed from $\oA$ in a very simple way.

In \textbf{Section \ref{SectionGS}}, we apply the theory to prove that the Gerstenhaber-Schack diagram cohomology is isomorphic to the operadic cohomology.
We begin by explaining how a diagram of associative algebras is described by an operad $\oA$.
We also make the associated module $\oMDA$ explicit.
Then we recall the Gerstenhaber-Schack cohomology and obtain a candidate for a resolution of $\oMDA$.
In Section \ref{SectionResolution}, we verify that the candidate is a valid resolution.
This computation is complicated, but still demonstrates the technical advantage of passing to the modules.

\section{Basics} \label{SectionBasics}

Fix the following symbols:
\begin{itemize}
\item $C$ is a set of colours.
\item $k$ is a field of characteristics $0$.
\item $\N_0$ is the set of natural numbers including $0$.
\end{itemize}

\noindent We will also use the following notations and conventions:
\begin{itemize}
\item Vector spaces over $k$ are called \emph{$k$-modules}, chain complexes of vector spaces over $k$ with differential of degree $-1$ are called \emph{dg-$k$-modules} and morphisms of chain complexes are called just \emph{maps}.
Chain complexes are assumed non-negatively graded unless stated otherwise.
\item $\dg{x}$ is the degree of an element $x$ of a dg-$k$-module.
\item $H_*(A)$ is homology of the object $A$, whatever $A$ is.
\item $\kspan{S}$ is the $k$-linear span of the set $S$.
\item $\ar{v}$ is arity of the object $v$, whatever $v$ is.
\item \emph{Quism} is a map $f$ of dg-$k$-modules such that the induced map $H_*(f)$ on homology is an isomorphism.
\end{itemize}

\begin{definition}
A \textbf{dg-\coll{C}} $X$ is a set $$\set{X\binom{c}{c_1,\ldots,c_n}} {n\in\N_0, c,c_1,\ldots,c_n\in C}$$ of dg-$k$-modules.
We call $c$ the \emph{output} colour of elements of $X\binom{c}{c_1,\ldots,c_n}$, $c_1,\ldots,c_n$ are the \emph{input} colours, $n$ is the \emph{arity}.
We also admit $n=0$.

When the above dg-$k$-modules have zero differentials, we talk just about graded \coll{C}.
If moreover no grading is given, we talk just about \coll{C}.
All notions that follow have similar analogues.
If the context is clear, we might omit the prefixes dg-$C$ completely.

A \textbf{dg-\colop{C}} $\oA$ is a dg-\coll{C} $\oA$ together with a set of of dg-$k$-module maps
$$\oo_i:\oA\binom{c}{c_1,\ldots,c_k}\ot\oA\binom{c_i}{d_1,\cdots,d_l}\to\oA\binom{c}{c_1,\cdots,c_{i-1},d_1,\cdots,d_l,c_{i+1},\cdots c_k},$$
called \emph{operadic compositions}, one for each choice of $k,l\in\N_0$, $1\leq i\leq k$ and $c,c_1,\ldots,c_k,d_1,\ldots,d_l \in C$, and a set of units
$$e:k\to\oA\binom{c}{c},$$
one for each $c\in C$.
These maps satisfy the usual associativity and unit axioms, e.g. \cite{MSS}.

The initial dg-\colop{C} is denoted $I$.

Equivalently, dg-\colop{C} is a monoid in the monoidal category of dg-\coll{C}s with the \textbf{composition product} $\oo$:
\begin{gather*}
(\oA\oo\oB)\binom{c}{c_1,\ldots,c_n}:=\\
\bigoplus_{\substack{k\geq 0,\\ i_1,\ldots,i_k\geq 0,\\ d_1,\ldots,d_k\in C}}\oA\binom{c}{d_1,\ldots,d_k}\ot\oB\binom{d_1}{c_1,\ldots,c_{i_1}}\ot\cdots\ot\oB\binom{d_k}{c_{i_1+\cdots i_{k-1}+1},\ldots,c_n}.
\end{gather*}

In contrast to the uncoloured operads, the composition is defined only for the ``correct'' colours and there is one unit in every colour, i.e. $I\binom{c}{c}=k$ for every $c\in C$.
Hence we usually talk about the unit\ul{s}.
We denote by $\id_c$ the image of $1\in k=I\binom{c}{c}$, hence $\Im e=\bigoplus_{c\in C}\kspan{\id_c}$.
The notation $\id_c$ for units coincides with the notation for identity morphisms.
The right meaning will always be clear from the context.

For a dg-\colop{C} $\oA$, we can consider its homology $H_*(\oA)$.
The operadic composition descends to $H_*(\oA)$.
Obviously, the units $\id_c$ are concentrated in degree $0$ and by our convention on non-negativity of the grading, $\id_c$ defines a homology class $[\id_c]$.
It is a unit in $H^*(\oA)$.
It can happen that $[\id_c]=0$ in which case it is easily seen that $H_*(\oA)\binom{c_0}{c_1,\ldots,c_n}=0$ whenever any of $c_i$'s equals $c$.
If all $[\id_c]$'s are nonzero, then $H_*(\oA)$ is a graded \colop{C}.

Let $\oM_1$ and $\oM_2$ be two dg-\coll{C}s.
Then \textbf{dg-\coll{C} morphism $f$} is a set of dg maps
$$f\binom{c}{c_1,\ldots,c_n} : \oM_1\binom{c}{c_1,\ldots,c_n}\to\oM_2\binom{c}{c_1,\ldots,c_n},$$
one for each $n\in\N_0, c,c_1,\ldots,c_n\in C$.

The dg-\coll{C} morphisms are composed ``colourwise'' in the obvious way.

A \textbf{dg-\colop{C} morphism} is a dg-\coll{C} morphisms preserving the operadic compositions and units.
\end{definition}

Recall that given a dg-$k$-module $A$, the \emph{endomorphism operad} $\oEnd_A$ is equipped with the differential
$$\dd_{\oEnd_A} f:=\dd_A f - (-1)^\dg{f} f\dd_{\otexp A n}$$
for $f\in\oEnd_A(n)$ homogeneous.
Let there be a decomposition $$(A,\dd_A)=\bigoplus_{c\in C}(A_c,\dd_{A_c}).$$
Then the endomorphism operad is naturally a dg-\colop{C} via
$$\oEnd_A\binom{c}{c_1,\ldots,c_n}:=\Hom_k(A_{c_1}\op\cdots\op A_{c_n},A_c).$$

An \textbf{algebra over a dg-\colop{C}} $\oA$ is a dg-\colop{C} morphism $$(\oA,\dd_{\oA})\to (\oEnd_A,\dd_{\oEnd_A}).$$

\subsection{Operadic modules} \label{OpMod}

\begin{definition}
Let $\oA=(\oA,\dd_\oA)$ be a dg-\colop{C}.
An (operadic) \textbf{dg-$\oA$-module} $\oM$ is a dg-\coll{C} 
$$\set{\oM\binom{c}{c_1,\ldots,c_n}}{n\in\N_0, c,c_1,\ldots,c_n\in C}$$
with structure maps
\begin{align*}
\oo^{L}_i &: \oA\binom{c}{c_1,\cdots,c_k} \ot \oM\binom{c_i}{d_1,\cdots,d_l} \to \oM\binom{c}{c_1,\cdots,c_{i-1},d_1,\cdots,d_l,c_{i+1},\cdots c_k},  \\
\oo^{R}_i &: \oM\binom{c}{c_1,\cdots,c_k} \ot \oA\binom{c_i}{d_1,\cdots,d_l} \to \oM\binom{c}{c_1,\cdots,c_{i-1},d_1,\cdots,d_l,c_{i+1},\cdots c_k}, \\
\end{align*}
one for each choice of $c,c_1,\cdots,d_1,\cdots \in C$ and $1\leq i \leq k$.
These structure maps are required to satisfy the expected axioms:
\begin{align*}
(\alpha_1 \oo_j \alpha_2) \oo_i^L m &= \left\{ \begin{array}{lcl}
																			 (-1)^{\dg{\alpha_2}\dg{m}}(\alpha_1\oo_i^L m)\oo_{j+\ar{m}-1}^R\alpha_2 & \ldots & i<j \\
																			 \alpha_1\oo_j^L(\alpha_2\oo_{i-j+1}^L m) & \ldots & j\leq i \leq j+\ar{\alpha_2}-1 \\
																			 (-1)^{\dg{\alpha_2}\dg{m}}(\alpha_1\oo_{i-\ar{\alpha_2}+1}^L)\oo_j^R\alpha_2 &  \ldots & i\geq j+\ar{\alpha_2},
																			 \end{array} \right. \\
m \oo_i^R (\alpha_1 \oo_j \alpha_2) &= (m \oo_i^R \alpha_1) \oo_{j+i-1}^R \alpha_2, \\
(\alpha_1 \oo_i^L m) \oo_j^R \alpha_2 &= \left\{ \begin{array}{lcl}
																				 (-1)^{\dg{\alpha_2}\dg{m}}(\alpha_1\oo_j\alpha_2)\oo_{i+\ar{\alpha_2}-1}^L m & \cdots & j<i \\
																				 \alpha_1\oo_i^L (m\oo_{j-i+1}^R \alpha_2) & \cdots & i\leq j\leq i+\ar{m}-1 \\
																				 (-1)^{\dg{\alpha_2}\dg{m}}(\alpha_1\oo_{j-\ar{m}+1}\alpha_2)\oo_i^L m & \cdots & j\geq i+\ar{m} \\
																				 \end{array} \right.
\end{align*}
and $$\id_c\oo_1 m=m=m\oo_i \id_{c_i}\quad\ldots\quad 1\leq i\leq\ar{m}$$
for $\alpha_1,\alpha_2 \in \oA$ and $m\in\oM\binom{c}{c_1,\cdots,{c_\ar{m}}}$ in the correct colours.
We usually omit the upper indices $L,R$, writing only $\oo_i$ for all the operations.

A \textbf{morphism of dg-$\oA$-modules} $\oM_1,\oM_2$ is a dg-\coll{C} morphism $\oM_1\xrightarrow{f}\oM_2$ satisfying
\begin{align*}
f(a \oo_i^L m) &= a \oo_i^L f(m), \\
f(m \oo_i^R a) &= f(m) \oo_i^R a.
\end{align*}

We expand the definition of dg-$\oA$-module.
Recall that each $\oM\binom{c}{c_1,\ldots,c_n}$ is a a dg-$k$-module.
The differentials of these dg-$k$-modules define a dg-\coll{C} morphism $\dd_\oM:\oM\to\oM$ of degree $-1$ satisfying $\dd^2_\oM=0$. 
The structure maps $\oo_i^L$ and $\oo_i^R$ commute with the differentials on the tensor products.
Hence $\dd_\oM:\oM\to\oM$ is a derivation in the following sense:
\begin{align*}
\dd_\oM(a\oo_i m) &= \dd_\oA a\oo_i m + (-1)^{\dg{a}} a\oo_i \dd_\oM m, \\
\dd_\oM(m\oo_i a) &= \dd_\oM m\oo_i a + (-1)^{\dg{m}} m\oo_i \dd_\oA a.
\end{align*}

As in the case of modules over a ring, $\oA$-modules form an abelian category.
We have ``colourwise'' kernels, cokernels, submodules etc.
There is a free $\oA$-module generated by a \coll{C} $M$, denoted
$$\FrM{\oA}{M}$$
and satisfying the usual universal property.
As an example of an explicit description of $\FrM{\oA}{M}$, let $\oA := \Fr{}{M_1}$ be a free \colop{C} generated by a \coll{C} $M_1$.
Then $\FrM{\oA}{M_2}$ is spanned by all planar trees, whose exactly one vertex is decorated by an element of $M_2$ and all the other vertices are decorated by elements of $M_1$ such that the colours are respected in the obvious sense.

We warn the reader that the notion of operadic module varies in the literature.
For example the monograph \cite{Modules} uses a different definition.
\bigskip

While dealing with $\oA$-modules, it is useful to introduce the following \textbf{infinitesimal composition product}  $\oA\oo'(\oB,\oC)$ of \coll{C}s $\oA,\oB,\oC$:
\begin{gather}
(\oA\oo'(\oB,\oC))\binom{c}{c_1,\ldots,c_n}:= \label{ICP} \\
\bigoplus_{\substack{k\geq 0,l>0,\\ 0\leq i_1\leq\ldots\leq i_{k-1}\leq n,\\ d_1,\ldots,d_k\in C}} \hspace{-0.75cm} \oA\binom{c}{d_1,\!...,d_k} \hspace{-0.15cm} \ot \hspace{-0.1cm} \oB\binom{d_1}{c_1,\!...,c_{i_1}} \hspace{-0.15cm} \ot \hspace{-0.1cm} \cdots \hspace{-0.1cm} \ot \hspace{-0.05cm} \underbrace{ \oC\binom{d_l}{c_{i_{l-1}+1},\!...,c_{i_{l}}} }_{l^{th}\mathrm{-position}} \hspace{-0.1cm} \ot \hspace{-0.1cm} \cdots \hspace{-0.1cm} \ot \hspace{-0.1cm} \oB\binom{d_k}{c_{i_{k-1}+1},\!...,c_n}. \nonumber
\end{gather}
See also \cite{LV}.
We denote by $$\oA\oo'_l(\oB,\oC)$$ the projection of $\oA\oo'(\oB,\oC)$ onto the component with fixed $l$.

For the free module, we have the following description using the infinitesimal composition product:
$$\FrM{\oA}{M}\cong\oA\oo'(I,M\oo\oA).$$
\end{definition}

\subsection{Tree composition} \label{TreeComposition}

An (unoriented) \emph{graph} (without loops) is a set $V$ of vertices, a set $H_v$ of half edges for every $v\in V$ and a set $E$ of (distinct) unordered pairs (called \emph{edges}) of distinct elements of $V$.
If $e:=(v,w)\in E$, we say that the vertices $v,w$ are adjacent to the edge $e$ and the edge $e$ is adjacent to the vertices $v,w$.
Denote $E_v$ the set of all edges adjacent to $v$.
Similarly, for $h\in H_v$, we say that the vertex $v$ is adjacent to the half edge $h$ and vice versa.
A path connecting vertices $v,w$ is a sequence $(v,v_1),(v_1,v_2),\ldots,(v_n,w)$ of \emph{distinct} edges.
A \emph{tree} is a graph such that for every two vertices $v,w$ there is a path connecting them iff $v\neq w$.
A \emph{rooted} tree is a tree with a chosen half edge, called \emph{root}.
The \emph{root vertex} is the unique vertex adjacent to the root.
The half edges other than the root are called \emph{leaves}.
For every vertex $v$ except for the root vertex, there is a unique edge $e_v\in E_v$ contained in the unique path connecting $v$ to the root vertex.
The edge $e_v$ is called \emph{output} and the other edges and half edges adjacent to $v$ are called \emph{legs} or \emph{inputs} of $v$.
The root is, by definition, the output of the root vertex.
The number of legs of $v$ is called arity of $v$ and is denoted $\ar{v}$.
Notice we also admit vertices with no legs, i.e. vertices of arity $0$.
A \emph{planar} tree is a tree with a given ordering of the set $H_v \sqcup E_v - \{e_v\}$ for each $v\in V$ (the notation $\sqcup$ stands for the disjoint union).
The planarity induces an ordering on the set of all leaves, e.g. by numbering them $1,2,\ldots$.
For example,
$$
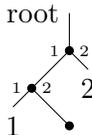
\begin{tikzpicture}[scale=0.5]
\draw (0,1) node[left]{root} -- (0,0) node[left]{$\scriptscriptstyle{1}$} node[right]{$\scriptscriptstyle{2}$} -- (-1,-1) node[left]{$\scriptscriptstyle{1}$} node[right]{$\scriptscriptstyle{2}$} -- (-1.5,-1.5) node[below]{$1$};
\draw (0,0) -- (0.5,-0.5) node[below]{$2$};
\draw (-1,-1) -- (0,-2);
\filldraw (0,0) circle (3pt);
\filldraw (-1,-1) circle (3pt);
\filldraw (0,-2) circle (3pt);
\end{tikzpicture}
$$
is a planar tree with $3$ vertices, $3$ half edges and $2$ edges.
We use the convention that the topmost half edge is always the root.
Then there are $2$ leaves.
The planar ordering of legs of all vertices is denoted by small numbers and the induced ordering of leaves is denoted by big numbers.

Let $T_1,T_2$ be two planar rooted trees, let $T_1$ have $n$ leaves and for $j=1,2$ let $V_j$ resp. $H_{j,v}$ resp. $E_j$ denote the set of vertices resp. half edges resp. edges of $T_j$.
For $1\leq i\leq n$ we have the grafting operation $\oo_i$ producing a planar rooted tree $T_1\oo_i T_2$ defined as follows: first denote $l$ the $i^{th}$ leg of $T_1$ and denote $v_1$ the vertex adjacent to $l$, denote $r$ the root of $T_2$ and denote $v_2$ the root vertex of $T_2$.
Then the set of vertices of $T_1\oo_i T_2$ is $V_1\sqcup V_2$, the set $H_v$ of half edges is
$$H_v=\left\{ \begin{array}{lcl} H_{1,v} & \ldots & v\in V_1 - \{v_1\}  \\ H_{1,v_1}-\{l\} & \ldots & v=v_1 \\ H_{2,v} & \ldots & v\in V_2-\{v_2\} \\ H_{2,v_2}-\{r\} & \ldots & v=v_2 \end{array} \right. ,$$
and finally the set of edges is $E_1\sqcup E_2\sqcup\{(v_1,v_2)\}$.
The planar structure is inherited in the obvious way.
For example,
$$
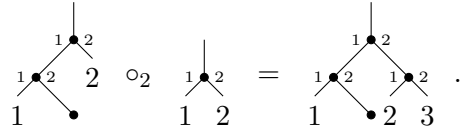
\raisebox{-0.7cm}{
\begin{tikzpicture}[scale=0.5]
\draw (0,1) -- (0,0) node[left]{$\scriptscriptstyle{1}$} node[right]{$\scriptscriptstyle{2}$} -- (-1,-1) node[left]{$\scriptscriptstyle{1}$} node[right]{$\scriptscriptstyle{2}$} -- (-1.5,-1.5) node[below]{$1$};
\draw (0,0) -- (0.5,-0.5) node[below]{$2$};
\draw (-1,-1) -- (0,-2);
\filldraw (0,0) circle (3pt);
\filldraw (-1,-1) circle (3pt);
\filldraw (0,-2) circle (3pt);
\end{tikzpicture}
}
\oo_2
\raisebox{-0.7cm}{
\begin{tikzpicture}[scale=0.5]
\draw (0,1) -- (0,0) node[left]{$\scriptscriptstyle{1}$} node[right]{$\scriptscriptstyle{2}$} -- (-0.5,-0.5) node[below]{$1$};
\draw (0,0) -- (0.5,-0.5) node[below]{$2$};
\filldraw (0,0) circle (3pt);
\end{tikzpicture}
}
=
\raisebox{-0.7cm}{
\begin{tikzpicture}[scale=0.5]
\draw (0,1) -- (0,0) node[left]{$\scriptscriptstyle{1}$} node[right]{$\scriptscriptstyle{2}$} -- (-1,-1) node[left]{$\scriptscriptstyle{1}$} node[right]{$\scriptscriptstyle{2}$} -- (-1.5,-1.5) node[below]{$1$};
\draw (0,0) -- (1,-1) node[left]{$\scriptscriptstyle{1}$} node[right]{$\scriptscriptstyle{2}$} -- (1.5,-1.5) node[below]{$3$};
\draw (1,-1) -- (0.5,-1.5) node[below]{$2$} ;
\draw (-1,-1) -- (0,-2);
\filldraw (0,0) circle (3pt);
\filldraw (-1,-1) circle (3pt);
\filldraw (0,-2) circle (3pt);
\filldraw (1,-1) circle (3pt);
\end{tikzpicture}
}.
$$

From this point on, \emph{tree} will always mean a planar rooted tree.
Such trees can be used to encode compositions of elements of an operad including those of arity $0$.

Let $T$ be a tree with $n$ vertices $v_1,\ldots,v_n$.
Suppose moreover that the vertices of $T$ are ordered, i.e. there is a bijection $b:\{v_1,\ldots,v_n\} \to \{1,\ldots,n\}$.
We denote such a tree with ordered vertices by $T_b$.

Now we explain how the bijection $b$ induces a structure of tree with levels on $T_b$ such that each vertex is on a different level.
Intuitively, $b$ encodes \emph{in what order} are elements of an operad composed.
We formalize this as follows:

Let $p_1,\ldots,p_n \in \oP$ be elements of a dg-\colop{C} $\oP$ such that if two vertices $v_i,v_j$ are adjacent to a common edge $e$, which is simultaneously the $l^{th}$ leg of $v_i$ and the output of $v_j$, then the $l^{th}$ input colour of $p_i$ equals the output colour of $p_j$.
We say that $v_i$ is \emph{decorated} by $p_i$.
Define inductively:
Let $i$ be such that $v_i$ is the root vertex.
Define
\begin{align*}
T^1 &:= v_i, \\
T^1_b(p_1,\ldots,p_n) &:= p_i.
\end{align*}
Here we are identifying $v_i$ with the corresponding corolla.
Assume a subtree $T^{k-1}$ of $T$ and $T^{k-1}_b(p_1,\ldots,p_n) \in \oP$ are already defined.
Consider the set $J$ of all $j$'s such that $v_j \not\in T^{k-1}$ and there is an edge $e$ between $v_j$ and some vertex $v$ in $T^{k-1}$.
Let $i\in J$ be such that $b(v_i)=\min\{b(v_j):j\in J\}$.
Let $l$ be the number of the leg $e$ of vertex $v$ in the planar ordering of $T$ and define
\begin{align*}
T^k &:= T^{k-1} \circ_l v_i, \\
T^k_b(p_1,\ldots,p_n) &:= T^{k-1}_b(p_1,\ldots,p_n) \circ_l p_i.
\end{align*}
In the upper equation, we are using the operation $\circ_l$ of grafting of trees.
Finally $$T_b(p_1,\ldots,p_n) := T^n_b(p_1,\ldots,p_n).$$

$T_b(p_1,\ldots,p_n)$ is called \textbf{tree composition} of $p_1,\ldots,p_n$ along $T_b$.

If $T$ and $p_i$'s are fixed, changing $b$ may change the sign of $T_b(p_1,\ldots,p_n)$.
Observe that if $\oP$ is concentrated in even degrees (in particular $0$) then the sign doesn't change.
If $b$ is understood and fixed, we usually omit it.

For example, let
$$T:=
\raisebox{-0.5cm}
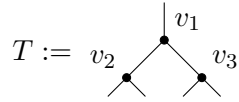
{\begin{tikzpicture}[scale=0.5]
\draw (0,1) -- (0,0) node[above right]{$v_1$} -- (-1,-1) node[above left]{$v_2$} -- (-1.5,-1.5);
\draw (-1,-1) -- (-0.5,-1.5);
\draw (0,0) -- (1,-1) node[above right]{$v_3$} -- (1.5,-1.5);
\draw (1,-1) -- (0.5,-1.5);
\filldraw (0,0) circle (3pt);
\filldraw (-1,-1) circle (3pt);
\filldraw (1,-1) circle (3pt);
\end{tikzpicture}
}
$$
For $i=1,2,3$, let $p_i$ be an element of degree $1$ and arity $2$.
Let $b(v_1)=1$, $b(v_2)=2$, $b(v_3)=3$ and $b'(v_1)=1$, $b'(v_2)=3$, $b(v_3)=2$.
Then $$T_{b}(p_1,p_2,p_3)=(p_1\oo_1 p_2)\oo_3 p_3\quad\mbox{and}\quad T_{b'}(p_1,p_2,p_3)=(p_1\oo_2 p_3)\oo_1 p_2$$
and by the associativity axiom $$T_{b}(p_1,p_2,p_3)=-T_{b'}(p_1,p_2,p_3).$$
\bigskip

A useful observation is that we can always reindex $p_i$'s so that
\begin{gather} \label{EasyTreeComp}
T_b(p_1,\ldots,p_n) = (\cdots((p_1\oo_{i_1}p_2)\oo_{i_2}p_3)\cdots\oo_{i_{n-1}}p_n)
\end{gather}
for some $i_1,i_2,\ldots,i_{n-1}$.

Tree compositions are a convenient notation for dealing with operadic derivations.

\subsection{Free product of operads} \label{SectionFreeProduct}

\begin{definition}
\textbf{Free product} $\oA\fpr\oB$ of dg-\colop{C}s $\oA,\oB$ is the coproduct $\oA\coprod\oB$ in the category of dg-\colop{C}s.
\end{definition}

Let $A,B$ be dg-$k$-modules.
The usual K\"unneth formula states that the map
\begin{align}
H_*(A)\ot H_*(B) &\xrightarrow{\iota} H_*(A\ot B) \label{UsualKunneth} \\
[a]\ot[b]&\mapsto [a\ot b] \nonumber
\end{align}
is a natural isomorphism of dg-$k$-modules, where $[\ ]$ denotes a homology class.
Our aim here is to prove an analogue of the K\"unneth formula for the free product of operads, that is 
\begin{gather} \label{KunnethFreeProd}
H_*(\oA)\fpr H_*(\oB)\cong H_*(\oA\fpr\oB)
\end{gather}
naturally as \colop{C}s.

First we describe $\oA\fpr\oB$ more explicitly.
Intuitively, $\oA\fpr\oB$ is spanned by trees whose vertices are decorated by elements of $\oA$ or $\oB$ such that no two vertices adjacent to a common edge are both decorated by $\oA$ or both by $\oB$.
Unfortunately, this is not quite true - there are problems with units of the operads.

Recall a dg-\colop{C} $\oP$ is called \emph{augmented} iff there is a dg-\colop{C} morphism $\oP\xrightarrow{a}I$ inverting the unit of $\oP$ on the left, i.e. the composition $I\xrightarrow{e}\oP\xrightarrow{a}I$ is $\id_I$.
The kernel of $a$ is denoted by $\ol{\oP}$ and usually called augmentation ideal.

If $\oA,\oB$ are augmented, we let the vertices be decorated by the augmentation ideals $\ol{\oA},\ol{\oB}$ instead of $\oA,\oB$ and the above description of $\oA\fpr\oB$ works well.
In fact, this has been already treated in \cite{HAAHA}.

However we will work without the augmentation assumption.
Choose a sub-\coll{C} $\ol{\oA}$ of $\oA$ such that
\begin{gather}
\ol{\oA}\op\id_{\oA} = \oA, \label{Condition1} \\
\Im\dd_\oA\subset\ol{\oA}, \label{Condition2}
\end{gather}
where $$\id_\oA:=\bigoplus_{c\in C}\kspan{\id_c}.$$
This is possible iff 
\begin{gather*} 
[\id_c]\neq 0 \mbox{ for all }c\in C,
\end{gather*}
that is if $H_*(\oA)$ is a graded \colop{C}.
This might not be the case generally as we have already seen at the beginning of Section \ref{SectionBasics} so let's assume it.
Choose $\ol{\oB}$ for $\oB$ similarly.

For given $\ol{\oA},\ol{\oB}$, a \textbf{free product tree} is a tree $T$ together with
$$c(v),c_1(v),c_2(v),\ldots,c_{\ar{v}}(v)\in C\mbox{ for each vertex }v$$
and a map 
\begin{gather} \label{oloPDefinition}
\ol{\oP}:\mbox{ vertices of } T\to\{\ol{\oA},\ol{\oB}\}
\end{gather}
such that if vertices $v_1,v_2$ are adjacent to a common edge, which is simultaneously the $l^{th}$ leg of $v_1$ and the output of $v_2$, then
$$c_l(v_1)=c(v_2) \quad\mbox{and}\quad \ol{\oP}(v_1)\neq \ol{\oP}(v_2).$$

Finally, the description of the free product is as follows:
\begin{gather} \label{FreeProductExplicitly}
\oA\fpr\oB := \bigoplus_{c\in C}\kspan{\id_c} \op \bigoplus_T\bigotimes_v \ol{\oP}(v)\binom{c(v)}{c_1(v),\ldots,c_{\ar{v}}(v)},
\end{gather}
where $T$ runs over all isomorphism classes of free product trees and $v$ runs over all vertices of $T$ and $\id_c$'s are of degree $0$.
If the vertices of $T$ are $v_1,\ldots,v_n$ then every element of $\bigotimes_{v} \ol{\oP}(v)\binom{c(v)}{c_1(v),\ldots,c_{\ar{v}}(v)}$ can be written as a tree composition $T(x_1,\ldots,x_n)$ where $x_i\in\ol{\oP}(v_i)$.
We say that $v_i$ is decorated by $x_i$.

The operadic composition $$T(x_1,\ldots,x_n)\oo_i T'(x'_1,\ldots,x'_m)$$ in $\oA\fpr\oB$ is defined in the obvious way by grafting $T$ and $T'$ (the result $T\oo_i T'$ of the grafting may not be a free product tree) and then (repeatedly) applying the following \emph{reducing operations}:
\begin{enumerate}
\item Suppose $w_1,w_2$ are vertices of $T\oo_i T'$ adjacent to a common edge $e$ which is simultaneously the $l^{th}$leg of $w_1$ and the output of $w_2$.
	Suppose moreover that $w_1$ is decorated by $p_1$ and $w_2$ by $p_2$.
	If both $p_1,p_2$ are elements of $\ol{\oA}$ or both of $\ol{\oB}$, then contract $e$ and decorate the resulting 			vertex by the composition of $p_1\oo_l p_2$.
\item If a vertex is decorated by a unit from $\id_\oA$ or $\id_\oB$ (this may happen since neither $\ol{\oA}$ nor $\ol{\oB}$ is generally closed under the composition!), omit it unless it is the only remaining vertex of the tree.
\end{enumerate}
After several applications of the above reducing operations, we obtain a free product tree or a tree with a single vertex decorated by a unit.

Obviously, $\id_c$'s are units for this composition.

The differential $\dd$ on $\oA\fpr\oB$ is determined by \eqref{FreeProductExplicitly} and the requirement that $\dd(\id_c)=0$ for every $c\in C$.
It has the derivation property and equals the differential on $\oA$ resp. $\oB$ upon the restriction on the corresponding sub-\colop{C} of $\oA\fpr\oB$.
Explicitly, for $T(x_1,\ldots,x_n)\in\oA\fpr\oB$ with $x_i\in\ol{\oA}\mbox{ or }\ol{\oB}$, assuming \eqref{EasyTreeComp}, we have
$$\dd(T(x_1,\ldots,x_n)) = \sum_{i=1}^n\epsilon_{i}T(x_1,\ldots,\dd(x_i),\ldots,x_n),$$
where $\epsilon_i:=(-1)^{\sum_{j=1}^{i-1}\dg{x_j}}$.

It is easily seen that the dg-\colop{C} $(\oA\fpr\oB,\dd)$ just described has the required universal property of the coproduct.

Now we are prepared to prove a version of \eqref{KunnethFreeProd} in a certain special case:

\begin{lemma} \label{LemmaKunneth}
Let $(\oA,\dd_\oA)\xrightarrow{\alpha}(\oA',\dd_{\oA'})$ and $(\oB,\dd_\oB)\xrightarrow{\beta}(\oB',\dd_{\oB'})$ be quisms of dg-\colop{C}s, that is we assume homology of $\oA,\oA',\oB,\oB'$ are graded \colop{C}s and $H_*(\alpha),H_*(\beta)$ are graded \colop{C} isomorphisms.
Then there are graded \colop{C} isomorphisms $\iota,\iota'$ such that the following diagram commutes:
\begin{diagram}
H_*(\oA)\fpr H_*(\oB) & \rTo^{\iota} & H_*(\oA\fpr\oB)  \\
\dTo<{H_*(\alpha)\fpr H_*(\beta)} & & \dTo>{H_*(\alpha\fpr\beta)} \\
H_*(\oA')\fpr H_*(\oB') & \rTo^{\iota'} & H_*(\oA'\fpr \oB')
\end{diagram}
\end{lemma}

\begin{proof}
Choose $\ol{\oA},\ol{\oB}$ so that \eqref{Condition1} and \eqref{Condition2} hold.
Now we want to choose $\ol{\oA}'\subset \oA'$ so that
\begin{gather}
\ol{\oA}'\op \id_{\oA'} = \oA', \nonumber \\
\alpha(\ol{\oA}) \subset \ol{\oA}' \label{ImageInOl}
\end{gather}
and choose $\ol{\oB}'\subset \oB'$ similarly.
To see that this is possible, we observe $\alpha(\ol{\oA})\cap\id_{\oA'}=0$:
If $\alpha(\ol{a})\in\id_{\oA'}$ for some $\ol{a}\in\ol{\oA}$, there is $u\in\id_\oA$ such that $\alpha(u)=\alpha(\ol{a})$, hence $\alpha(\ol{a}-u)=0$ and $\dd_{\oA}(\ol{a}-u)=0$ since both $\ol{a}$ and $u$ are of degree $0$.
Since $\alpha$ is a quism, $\ol{a}-u=\dd_{\oA}a$ for some $a\in\oA$ and by the property \eqref{Condition2} of $\ol{\oA}$ we have $\ol{a}-u\in\ol{\oA}$.
But this implies $u\in\ol{\oA}$, a contradiction.

Now use the explicit description \eqref{FreeProductExplicitly} of the free product $\oA\fpr\oB$ and the usual K\"unneth formula \eqref{UsualKunneth} to obtain an isomorphism
$$H_*(\oA)\fpr H_*(\oB) = \bigoplus_T\bigotimes_v H_*(\ol{\oP}(v)) \xrightarrow{\iota} H_*(\bigoplus_T\bigotimes_v \ol{\oP}(v)) = H_*(\oA\fpr\oB)$$
and similarly for $\oA',\oB'$.

Assume we are given a free product tree $T$ and its vertex $v$.
The tree $T$ comes equipped with $\ol{\oP}$ as in \eqref{oloPDefinition}.
Let $$\ol{\oP}'(v):=(\ol{\oP}(v))'=\left\{ \begin{array}{lcl} \ol{\oA}' & \mbox{ for } & \ol{\oP}(v)=\ol{\oA} \\ \ol{\oB}' & \mbox{ for } & \ol{\oP}(v)=\ol{\oB} \end{array} \right.$$
and define a map 
\begin{gather*}
\pi(v):\ol{\oP}(v)\to\ol{\oP}'(v), \\
\pi(v)=\left\{ \begin{array}{lcl} \alpha & \mbox{ for } & \ol{\oP}(v)=\ol{\oA} \\ \beta & \mbox{ for }  & \ol{\oP}(v)=\ol{\oB} \end{array} \right. .
\end{gather*}
This is justified by \eqref{ImageInOl}.
Then the following diagram
\begin{diagram}
\bigoplus_T\bigotimes_v H_*(\ol{\oP}(v)) & \rTo^{\iota} & H_*(\bigoplus_T\bigotimes_v \ol{\oP}(v))\\
\dTo<{\bigoplus_T\bigotimes_v H_*(\pi(v))} & & \dTo>{H_*(\bigoplus_T\bigotimes_v \pi(v))} \\
\bigoplus_T\bigotimes_v H_*(\ol{\oP}'(v)) & \rTo^{\iota'} & H_*(\bigoplus_T\bigotimes_v \ol{\oP}'(v))
\end{diagram}
commutes by the naturality of the usual K\"unneth formula.
The horizontal \coll{C} isomorphism $\iota$ is given in terms of tree compositions by the formula
$$\iota(T([x_1],[x_2],\ldots)) = [T(x_1,x_2,\ldots)],$$
where $x_1,x_2,\ldots\in\ol{\oA}$ or $\ol{\oB}$.
Now we verify that $\iota$ preserves the operadic composition:
$$\iota(T_x([x_1],\ldots)) \oo_i \iota(T_y([y_1],\ldots)) = \iota(T_x([x_1],\ldots)\oo_i T_y([y_1],\ldots)).$$
The left-hand side equals $[T_x(x_1,\ldots) \oo_i T_y(y_1,\ldots)]$, so we check
$$[T_x(x_1,\ldots) \oo_i T_y(y_1,\ldots)]=\iota(T_x([x_1],\ldots)\oo_i T_y([y_1],\ldots)).$$
We would like to perform the same reducing operations on $T_x(x_1,\ldots) \oo_i T_y(y_1,\ldots)$ and $T_x([x_1],\ldots)\oo_i T_y([y_1],\ldots)$ parallely.
For the first reducing operation, this is OK.
For the second one, if, say, $[x_1]\in\id_{H_*(\oA)}$, then $x_1=u+\dd_\oA a$ for some $u\in\id_\oA$ and $a\in\oA$.
Hence $T_x(x_1,\ldots)=T_x(u,\ldots)+T_x(\dd_\oA a,\ldots)$.
So we can go on with $T_x(u,\ldots) \oo_i T_y(y_1,\ldots)$ and $T_x([x_1],\ldots)\oo_i T_y([y_1],\ldots)$, ommiting the vertex $v_1$ decorated by $u$ resp. $[x_1]$, but we also have to apply the reducing operations to $T_x(\dd_\oA a,\ldots) \oo_i T_y(y_1,\ldots)$.
As it turns out, this tree composition is a boundary in $\oA\fpr\oB$.
We leave the details to the reader.
\end{proof}

\section{Operadic cohomology of algebras} \label{SectionOC}

\subsection{Reminder} \label{REMINDER}

Let $(\oR,\dd_\oR) \xrightarrow{\rho} (\oA,\dd_\oA)$ be dg-\colop{C} over $(\oA,\dd_\oA)$, i.e. $\rho$ is a dg-\colop{C} morphism.
Let $(\oM,\dd_\oM)$ be a dg-$\oA$-module.
Define a $k$-module $$\Der_\oA^n(\oR,\oM)$$
consisting of all \coll{C} morphisms $\theta:\oR\to\oM$ of degree $\dg{\theta}=n$ in all colours satisfying
$$ \theta(r_1 \circ_i r_2) = \theta(r_1)\oo_i^R \rho(r_2) + (-1)^{\dg{\theta}\dg{r_1}} \rho(r_1)\oo_i^L \theta(r_2) $$
for any $r_1,r_2\in \oR$ and any $1\leq i\leq \ar{r_1}$.
Denote
$$\Der_\oA (\oR,\oM) := \bigoplus_{n\in \Z} \Der_\oA^n (\oR,\oM).$$

For $\theta\in\Der_\oA (\oR,\oM)$ homogeneous, let
\begin{gather} \label{CotangDiff}
\codd\theta:=\theta\dd_\oR - (-1)^\dg{\theta}\dd_\oM\theta.
\end{gather}
Extending by linearity, the above formula defines a map $\delta$ from the $k$-module $\Der_\oA (\oR,\oM)$.

\begin{lemma}
$\codd$ maps derivations to derivations and $\codd^2=0$.
\end{lemma}

\begin{proof}
The degree of $\delta$ obviously equals $-1$ and
\begin{align*}
\codd^2\theta &= (\theta\dd_\oR - (-1)^\dg{\theta}\dd_\oM\theta)\dd_\oR -(-1)^{\dg{\delta\theta}}\dd_\oM(\theta\dd_\oR - (-1)^\dg{\theta}\dd_\oM\theta) = \\
&= \theta\dd_\oR^2 - (-1)^\dg{\theta}\dd_\oM\theta\dd_\oR -(-1)^{\dg\theta +1}\dd_\oM\theta\dd_\oR -(-1)^{\dg{\theta}+1+\dg\theta +1}\dd_\oM^2\theta = \\
&= 0.
\end{align*}
The following computation shows that $\codd$ maps derivations to derivations:
\begin{align*}
(\codd\theta)(r_1\oo_i r_2) &= \theta \left( \dd r_1\oo_i r_2 +(-1)^\dg{r_1}r_1 \oo_i \dd r_2 \right) + \\
&\phantom{=} -(-1)^\dg\theta \dd \left( \theta r_1 \oo_i \rho r_2 + (-1)^{\dg{\theta}\dg{r_1}} \rho r_1\oo_i\theta r_2 \right) = \\
&= \theta\dd r_1 \oo_i \rho r_2 +(-1)^{\dg{\theta}(\dg{r_1}+1)} \rho\dd r_1 \oo_i \theta r_2 + \\
&\phantom{=} +(-1)^{\dg{r_1}} \theta r_1 \oo_i \rho\dd r_2 + (-1)^{(\dg{\theta}+1)\dg{r_1}} \rho r_1 \oo_i \theta \dd r_2 + \\
&\phantom{=} -(-1)^\dg{\theta}\dd\theta r_1 \oo_i \rho r_2 - (-1)^{\dg{r_1}} \theta r_1\oo_i \dd\rho r_2 + \\
&\phantom{=} -(-1)^{\dg{\theta}(\dg{r_1}+1)} \dd\rho r_1\oo_i\theta r_2 - (-1)^{(\dg{\theta}+1)\dg{r_1}+\dg{\theta}}  \rho r_1\oo_i \dd\theta r_2 = \\
&= (\codd\theta)r_1\oo_i \rho r_2 + (-1)^{\dg{\codd\theta}\cdot\dg{r_1}} \rho r_1\oo_i (\codd\theta)r_2
\end{align*}
where we have ommited the subscripts of $\dd$.
\end{proof}

A particular example of this construction is
$$(\oR,\dd_\oR) \xrightarrow[\rho]{\sim} (\oA,\dd_\oA),$$
a cofibrant \cite{HAAHA} resolution of a dg-\colop{C} $\oA$, and
$$\oM:=(\oEnd_{\oA},\dd_{\oEnd_A}),$$
which is a dg-$\oA$-module via a dg-\colop{C} morphism
$$(\oA,\dd_\oA) \xrightarrow{\alg} (\oEnd_{A},\dd_{\oEnd_A})$$ determining an $\oA$-algebra structure on a dg-$k$-module $(A,\dd_A)=\bigoplus_{c\in C}(A_c,\dd_{A_c})$.

Let $\upar C$ denote the \emph{suspension} of a graded object $C$, that is $(\upar C)_n:=C_{n-1}$.
Analogously $\downar$ denote the \emph{desuspension}.

\begin{definition}
\begin{gather} \label{CotangDef}
\left( C^*(A,A),\codd \right) := \upar \left( \Der^{-*}((\oR,\dd_\oR), \oEnd_A),\codd \right)
\end{gather}
is called \textbf{operadic cotangent complex of the $\oA$-algebra $A$} and
$$H^*(A,A) := H^*(C^*(A,A),\codd)$$ is called \textbf{operadic cohomology of $\oA$-algebra $A$}.
\end{definition}

The change of grading $*\mapsto 1-*$ is purely conventional.
For example, if $\oA$ is the operad for associative algebras and $\oR$ is its minimal resolution, under our convention we recover the grading of the Hochschild complex for which the bilinear cochains are of degree $1$.

\subsection{Algebras with derivation} \label{AlgebrasWithDerivation}

Let $\oA$ be a dg-\colop{C}.
Consider a \coll{C} $\Phi := \kspan{\phi_c | c\in C}$, such that $\phi_c$ is of arity 1, degree $0$ and the input and output colours are both $c$.
Let $\mathfrak{D}$ be the ideal in $\oA \fpr \Fr{}{\Phi}$ generated by all elements
\begin{gather} \label{relator}
\phi_c \oo_1 \alpha - \sum_{i=1}^{n} \alpha\oo_i \phi_{c_i}
\end{gather}
for $n\in\N_0, c,c_1,\ldots,c_n\in C$ and $\alpha\in\oA\binom{c}{c_1,\ldots,c_n}$.
Denote $$\oDA := \left( \frac{\oA \fpr \Fr{}{\Phi}}{\mathfrak{D}} , \ddDA \right),$$
where $\ddDA$ is the derivation given by the formulas
$$\ddDA(a) := \ddA(a),\quad \ddDA(\phi_c) := 0$$
for $a\in\oA$ and $c\in C$.

An algebra over $\oDA$ is a pair $(A,\phi)$, where $A=\bigoplus_{c\in C}A_c$ is an algebra over $\oA$ and $\phi$ is a derivation of $A$ in the following sense:
$\phi$ is a collection of degree $0$ dg-maps $\phi_c:A_c\to A_c$ such that
$$\phi_c(\alpha(a_1,\ldots,a_n)) = \sum_{i=1}^n \alpha(a_1,\ldots,\phi_{c_i}(a_i),\ldots,a_n)$$
for $\alpha\in\oA\binom{c}{c_1,\ldots,c_n}$ and $a_j\in A_{c_j}$, $1\leq j\leq n$.

Given a \emph{free} resolution
\begin{gather} \label{resolutionR}
\oR := (\Fr{}{X},\ddR) \xrightarrow{\rR} (\oA,\ddA),
\end{gather}
where $X$ is a dg-\coll{C}, it is surprisingly easy to explicitly construct a free resolution of $(\oDA,\ddDA)$.
Consider the free graded \colop{C} 
$$\oDR := \Fr{}{X\op\Phi\op\ul{X}},$$
where $\underline{X}:=\upar X$.
We denote by $\ul{x}$ the element $\upar x\in\ul{X}$ corresponding to $x\in X$.
To describe the differential, let $s:\Fr{}{X} \to \oDR$ be a degree $+1$ derivation determined by
$$s(x) := \ul{x} \mbox{ for }x\in X.$$
Then define a degree $-1$ derivation $\ddDR : \oDR \to \oDR$ by
\begin{align}
&\ddDR(x) := \ddR(x), \nonumber \\
&\ddDR(\phi_c) := 0, \label{DiffOnoDR} \\ 
&\ddDR(\ul{x}) := \phi_c \oo_1 x - \sum_{i=0}^n x\oo_i \phi_{c_i} - s(\ddR x). \nonumber
\end{align}

\begin{convention}
From now on we will assume
$$n\in \N_0,\ c,c_1,\ldots,c_n\in C,\ x\in X\binom{c}{c_1,\ldots,c_n}$$
whenever any of these symbols appears.
We will usually omit the lower indices $c$ and $c_i$'s for $\phi$.
\end{convention}

\begin{lemma}
$\ddDR^2 = 0$.
\end{lemma}

\begin{proof}
Using the tree compositions of Section \ref{TreeComposition}, let $\ddR(x)=\sum_i T_i(x_{i1},\cdots,x_{in_i})$.
\begin{align*}
\ddDR^2(x) &= \ddDR \left( \phi\oo_1 x - \sum_{j=1}^n x \oo_j\phi - s(\ddR(x)) \right) = \\
&= \phi\oo_1 \ddDR (x) - \sum_{j=1}^n \ddDR (x) \oo_j\phi\ + \\
&\phantom{=} - \ddDR \left( \sum_i \sum_{j=1}^{n_i} \epsilon_{ij} T_i(x_{i1},\ldots,\ul{x_{ij}},\ldots,x_{in_i}) \right) 
\end{align*}
If we assume \eqref{EasyTreeComp}, then $\epsilon_{ij}=(-1)^{\sum_{l=1}^{j-1}\dg{x_{il}}}$.
The last application of $\ddDR$ on the double sum can be rewritten as
\begin{align*}
&{} \sum_i \sum_{j=1}^{n_i} \sum_{\substack{1\leq k\leq n_i,\\ k\neq j}} \tilde{\epsilon}_{ijk} T_i(x_{i1},\ldots,\ddR(x_{ik}),\ldots,\ul{x_{ij}},\ldots,x_{in_i})\ + \\
&{} + \sum_i \sum_{j=1}^{n_i} T_i(x_{i1},\ldots,\phi\oo_1 x_{ij},\ldots,x_{in_i})\ + \\
&{} - \sum_i \sum_{j=1}^{n_i} \sum_{k=1}^{\ar{x_{ij}}} T_i(x_{i1},\ldots, x_{ij}\oo_k \phi,\ldots,x_{in_i})\ + \\
&{} - \sum_i \sum_{j=1}^{n_i} T_i(x_{i1},\ldots,s(\ddR(x_{ij})),\ldots,x_{in_i}),
\end{align*}
where $\tilde{\epsilon}_{ijk}=\epsilon_{ij}\epsilon_{ik}$ if $k<j$ and $\tilde{\epsilon}_{ijk}=-\epsilon_{ij}\epsilon_{ik}$ if $k>j$.
The second and third lines sum to 
\begin{gather*}
\sum_i \phi\oo_1 T_i(x_{i1},\ldots,x_{in_i}) - \sum_i \sum_{j=1}^\ar{T_i} T_i(x_{i1},\ldots,x_{in_i})\oo_j \phi = \\
= \phi\oo_1 \ddDR (x) - \sum_{j=1}^n \ddDR (x) \oo_j\phi,
\end{gather*}
while the first and last rows sum to 
$$-s \ddDR \left( \sum_i T_i(x_{i1},\cdots,x_{in_i}) \right) = -s\ddDR^2 = 0$$
and this concludes the computation.
\end{proof}

From now on, we will refer by $\oDR$ also to the \emph{dg}-\colop{C} $(\oDR,\ddDR)$.
Define a \colop{C} morphism $\rDR : \oDR \to \oDA$ by
\begin{align*}
&\rDR (x) := \rR(x), \\
&\rDR (\phi_c) := \phi_c, \\
&\rDR (\ul{x}) := 0. \\
\end{align*}

\begin{theorem} \label{ResolutionOfDerivation}
$\rDR$ is a free resolution of $\oDA$.
\end{theorem}

\begin{example} \label{ExampleDR}
Let's see what we get for $\oA:=\oAss=\Fr{}{\mu}/(\mu\oo_1\mu-\mu\oo_2\mu)$ and its minimal resolution (see e.g. \cite{HDA}) $\oR:=\oAss_\infty=(\Fr{}{X},\ddR)\xrightarrow{\rR}(\oAss,0)$, where
$$X=\kspan{x^2,x^3,\ldots}$$
is the collection spanned by $x^n$ in arity $n$ and degree $\dg{x^n}=n-2$ and $\ddR$ is a derivation differential given by 
$$\ddR(x^n):=\sum_{i+j=n+1}\sum_{k=1}^i (-1)^{i+(k+1)(j+1)}x^i\oo_k x^j$$ 
and the quism $\rR:\oR=\oAss_\infty\to\oAss=\oA$ is given by 
$$\rR(x^2):=\mu, \quad \rR(x^n):=0\mbox{ for }n\geq 3.$$
Then the associated operad with derivation is
$$\oDA:=\frac{\oAss\fpr \Fr{}{\Phi}}{(\phi\oo\mu-\mu\oo_1\phi-\mu\oo_2\phi)},$$
where $\Phi:=\kspan{\phi}$ with $\phi$ a generator of arity $1$.
Its free resolution is $$\oDR:=(\Fr{}{X\op\Phi\op\ul{X}},\ddDR)\xrightarrow{\rDR}(\oDA,0),$$
where the differential $\ddDR$ is given by
\begin{align*}
&\ddDR(x) := \ddR(x),\\
&\ddDR(\phi) := 0, \\
&\ddDR(\ul{x}^n) := \phi\oo_1 x^n - \sum_{i=1}^n x^n\oo_i\phi - \sum_{i+j=n+1}\sum_{k=1}^i (-1)^{i+(k+1)(j+1)} (\ul{x}^i\oo_k x^j + (-1)^i x^i\oo_k\ul{x}^j)
\end{align*}
and the quism $\rDR$ by
$$\rDR(x):=\rR(x),\quad \rDR(\phi):=\phi,\quad \rDR(\ul{x})=0.$$
\end{example}

\begin{proof}[of Theorem \ref{ResolutionOfDerivation}]
Obviously $\rDR$ has degree $0$ and commutes with differentials because of the relations in $\oDA$.
Let's abbreviate $\ddDR =: \dd$.
First we want to use a spectral sequence to split $\dd$ such that $\dd^0$, the $0^{th}$ page part of $\dd$, is nontrivial only on the generators from $\ul{X}$.

Let's put an additional grading $\gr$ on the \coll{C} $X\op\Phi\op\ul{X}$ of generators:
$$\gr(x):=\dg{x},\quad \gr(\phi):=1,\quad\gr(\ul{x}):=\dg{\ul{x}}.$$
This induces a grading on $\oDR$ determined by the requirement that the composition is of $\gr$ degree $0$.
Let $$\mathfrak{F}_p:=\bigoplus_{i=0}^p\set{z\in\oDR}{\gr(z)=i}.$$
Obviously $\ddDR\mathfrak{F}_p\subset\mathfrak{F}_p$.
Consider the spectral sequence $E^*$ associated to the filtration
$$0\into\mathfrak{F_0}\into\mathfrak{F_1}\into\cdots$$
of $\oDR$.
On $\oDA$ we have the trivial filtration $$0\into\oDA$$ and the associated spectral sequence $E'^*$.

We will show that $\rDR$ induces quism $(E^1,\dd^1)\xrightarrow{\sim} (E'^1,\dd'^1)$.
Then we can use the comparison theorem since both filtrations are obviously bounded below and exhaustive (e.g. \cite{Weibel}, page 126, Theorem 5.2.12, and page 135, Theorem 5.5.1).

Then the $0^{th}$ page satisfies $E^0\cong\Fr{}{X\op\Phi\op\ul{X}}$ and is equipped with the derivation differential $\dd^0$:
$$\dd^0(x)=0=\dd^0(\phi_c),\quad \dd^0(\ul{x})=\phi_c\oo_1 x-\sum_{i=1}^\ar{x}x\oo_i \phi_{c_i}.$$
Denote by $\mathfrak{D}$ the ideal in $\Fr{}{X\op\Phi}$ generated by
\begin{gather} \label{DerivationRelation}
\phi_c\oo_1 x-\sum_{i=1}^\ar{x}x\oo_i \phi_{c_i}
\end{gather}
for all $x\in X\binom{c}{c_1,\ldots,c_\ar{x}}$ of arbitrary colours.

\begin{sublemma} \label{sublemma0}
$$H_*(E^0,\dd^0) \cong \frac{ \Fr{}{X\op\Phi} } { \mathfrak{D} }$$
\end{sublemma}

Once this sublemma is proved, $\dd^1$ on $E^1\cong H_*(E^0,\dd^0)\cong \frac{ \Fr{}{X\op\Phi} } { \mathfrak{D} }$ will be given by
\begin{gather} \label{CheckLater}
\dd^1(x)=\ddR(x),\quad \dd^1(\phi_c)=0.
\end{gather}
We immediately see that $E'^1\cong \oDA$ and it is equipped with the differential $\dd'^1=\ddDA$.
To see that $\rDR^1:E^1\to E'^1$ induced by $\rDR$ is a quism, observe that we can use the relations \eqref{relator} in $\oDA$ to "move all the $\phi$'s to the bottom of the tree compositions", hence, denoting $$\Phi':=\Fr{}{\Phi},$$
we have
$$\oDA \cong \oA\oo\Phi'.$$
The composition and the differential on $\oA\oo\Phi'$ are transferred along this isomorphism from $\oDA$.
Similarly,
\begin{gather} \label{XPhiQuot}
\frac{\Fr{}{X\op{\Phi}}}{\mathfrak{D}} \cong \Fr{}{X}\oo\Phi'.
\end{gather}
Under these quisms
\begin{gather} \label{CheckLater2}
\rDR^1 \mbox{ becomes } \rR\oo\id_{\Phi'}.
\end{gather}
It remains to use the usual K\"unneth formula \eqref{UsualKunneth} to finish the proof.
\bigskip

\begin{proof}[of Sublemma \ref{sublemma0}]
Denote $\phi_c^m := \phi_c\oo_1\cdots\oo_1\phi_c$ the $m$-fold composition of $\phi_c$.
Let $$\oDR^0 := \Fr{}{X\op\Phi}$$
and, for $n\geq 0$, let $\oDR^{n+1}\subset\oDR$ be spanned by elements
\begin{align*}
&\phi_c^m\oo_1 x \oo (x_1,\ldots,x_\ar{x}) \mbox{ and} \\
&\phi_c^m\oo_1 \ul{x} \oo (x_1,\ldots,x_\ar{x})
\end{align*}
for all $x\in X\binom{c}{c_1,\ldots,c_\ar{x}}, m\geq 0, x_i\in \oDR^n\binom{c_i}{\cdots}, 1\leq i\leq \ar{x}$.
In other words,
$$\oDR^{n+1} = \Phi'\oo(X\op\ul{X})\oo\oDR^n.$$
$\oDR^n$ is obviously closed under $\dd^0$ and
$$\oDR^0\into\oDR^1\into\cdots\to\colim_n \oDR^n\cong\oDR,$$
where the colimit is taken in the category of dg-\coll{C}s.

Before we go further, we must make a short notational digression.
Let $T(g_1,\ldots,g_m)$ be a tree composition with $g_i\in X\sqcup \Phi\sqcup \ul{X}$ for $1\leq i \leq m$.
Recall the tree $T$ has vertices $v_1,\ldots,v_m$ decorated by $g_1,\ldots,g_m$ (in that order).
We say that \textbf{$g_j$ is in depth $d$ in $T(g_1,\ldots,g_j,\ldots,g_m)$} iff the shortest path from $v_j$ to the root vertex passes through exactly $d$ vertices (including $v_j$ and the root vertex) decorated by elements of $X\sqcup\ul{X}$.

As an example, consider
$$
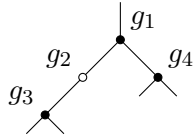
\begin{tikzpicture}[scale=0.5]
\draw (0,1) -- (0,0) node[above right]{$g_1$} -- (-1,-1) node[above left]{$g_2$} -- (-2,-2) node[above left]{$g_3$} -- (-2.5,-2.5);
\draw (-2,-2) -- (-1.5,-2.5);
\draw (0,0) -- (1,-1) node[above right]{$g_4$} -- (1.5,-1.5);
\draw (1,-1) -- (0.5,-1.5);
\filldraw (0,0) circle (3pt);
\filldraw[fill=white] (-1,-1) circle (3pt);
\filldraw (-2,-2) circle (3pt);
\filldraw (1,-1) circle (3pt);
\end{tikzpicture}
$$
If $g_1,g_3,g_4\in X$ and $g_2\in\Phi$, then $g_1,g_2$ are in depth $1$ and $g_3,g_4$ are in depth $2$.

Using the notion of depth, the definition of $\oDR^n$ can be rephrased as follows:
$\oDR^n$ is spanned by $T(g_1,\ldots,g_m)$ with $g_i\in X\sqcup\Phi\sqcup\ul{X},\ 1\leq i\leq m$, such that if $g_j\in\ul{X}$ for some $j$, then $g_j$ is in depth $\leq n$ in $T(g_1,\ldots,g_m)$.

Consider the quotient $\mathcal{Q}^n$ of $\Fr{}{X\op\Phi}$ by the ideal generated by elements
$$T(g_1,\ldots,g_{j-1},\phi_c\oo_1 x_j - \sum_{i=0}^\ar{x_j}x_j\oo_i \phi_{c_i},g_{j+1},\ldots,g_m)$$
for any tree $T$, any $g_1,\ldots,g_{j-1},g_{j+1},\ldots,g_m\in X\sqcup\Phi$ and any $x_j\in X$ in depth $\leq n$ in $T(g_1,\ldots,g_{j-1},x_j,g_{j+1},\ldots,g_m)$.
There are obvious projections
$$\Fr{}{X\op\Phi}=\mathcal{Q}^0 \onto \mathcal{Q}^1 \onto \cdots\to\colim_n \mathcal{Q}^n\cong \frac{\Fr{}{X\op\Phi}}{\mathfrak{D}}.$$
To see the last isomorphism, observe that we can use the relations defining $\mathcal{Q}^n$ to "move" the $\phi_c$'s in tree compositions so that they are all in depth $\geq n$ or in positions such that their inputs are leaves, then use \eqref{XPhiQuot}.

For example, consider the following computation in $\mathcal{Q}^2$, where the black vertices are decorated by $X$ and white vertices by $\Phi$:
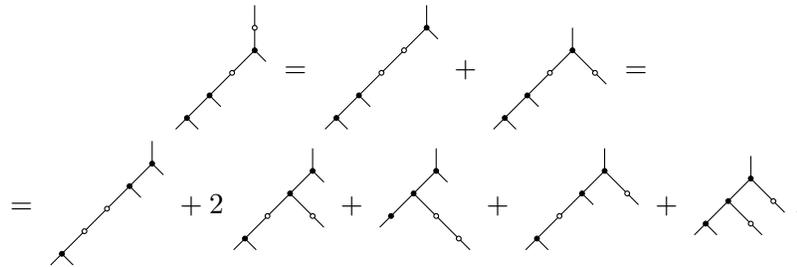
\begin{gather*}
\raisebox{-0.7cm}{
\begin{tikzpicture}[scale=0.3]
\draw (0,2) -- (0,0) -- (-1,-1) -- (-2,-2) -- (-3,-3) -- (-3.5,-3.5);
\draw (-2,-2) -- (-1.5,-2.5);
\draw (-3,-3) -- (-2.5,-3.5);
\draw (0,0) -- (0.5,-0.5);
\filldraw[fill=white] (0,1) circle (3pt);
\filldraw (0,0) circle (3pt);
\filldraw[fill=white] (-1,-1) circle (3pt);
\filldraw (-2,-2) circle (3pt);
\filldraw (-3,-3) circle (3pt);
\end{tikzpicture}
}
=
\raisebox{-0.7cm}{
\begin{tikzpicture}[scale=0.3]
\draw (0,1) -- (0,0) -- (-1,-1) -- (-2,-2) -- (-3,-3) -- (-4,-4) -- (-4.5,-4.5);
\draw (-3,-3) -- (-2.5,-3.5);
\draw (-4,-4) -- (-3.5,-4.5);
\draw (0,0) -- (0.5,-0.5);
\filldraw (0,0) circle (3pt);
\filldraw[fill=white] (-1,-1) circle (3pt);
\filldraw[fill=white] (-2,-2) circle (3pt);
\filldraw (-3,-3) circle (3pt);
\filldraw (-4,-4) circle (3pt);
\end{tikzpicture}
}
+
\raisebox{-0.7cm}{
\begin{tikzpicture}[scale=0.3]
\draw (0,1) -- (0,0) -- (-1,-1) -- (-2,-2) -- (-3,-3) -- (-3.5,-3.5);
\draw (-2,-2) -- (-1.5,-2.5);
\draw (-3,-3) -- (-2.5,-3.5);
\draw (0,0) -- (1,-1) -- (1.5,-1.5);
\filldraw (0,0) circle (3pt);
\filldraw[fill=white] (-1,-1) circle (3pt);
\filldraw (-2,-2) circle (3pt);
\filldraw (-3,-3) circle (3pt);
\filldraw[fill=white] (1,-1) circle (3pt);
\end{tikzpicture}
}
=\\
=
\raisebox{-0.7cm}{
\begin{tikzpicture}[scale=0.3]
\draw (0,1) -- (0,0) -- (-1,-1) -- (-2,-2) -- (-3,-3) -- (-4,-4) -- (-4.5,-4.5);
\draw (0,0) -- (0.5,-0.5);
\draw (-1,-1) -- (-0.5,-1.5);
\draw (-4,-4) -- (-3.5,-4.5);
\filldraw (0,0) circle (3pt);
\filldraw (-1,-1) circle (3pt);
\filldraw[fill=white] (-2,-2) circle (3pt);
\filldraw[fill=white] (-3,-3) circle (3pt);
\filldraw (-4,-4) circle (3pt);
\end{tikzpicture}
}
+2
\raisebox{-0.5cm}{
\begin{tikzpicture}[scale=0.3]
\draw (0,1) -- (0,0) -- (-1,-1) -- (-2,-2) -- (-3,-3) -- (-3.5,-3.5);
\draw (0,0) -- (0.5,-0.5);
\draw (-1,-1) -- (0,-2);
\draw (-3,-3) -- (-2.5,-3.5);
\draw (0,-2) -- (0.5,-2.5);
\filldraw (0,0) circle (3pt);
\filldraw (-1,-1) circle (3pt);
\filldraw[fill=white] (-2,-2) circle (3pt);
\filldraw (-3,-3) circle (3pt);
\filldraw[fill=white] (0,-2) circle (3pt);
\end{tikzpicture}
}
+
\raisebox{-0.5cm}{
\begin{tikzpicture}[scale=0.3]
\draw (0,1) -- (0,0) -- (-1,-1) -- (-2,-2) -- (-2.5,-2.5);
\draw (0,0) -- (0.5,-0.5);
\draw (-1,-1) -- (0,-2);
\draw (0,-2) -- (1,-3) -- (1.5,-3.5);
\filldraw (0,0) circle (3pt);
\filldraw (-1,-1) circle (3pt);
\filldraw (-2,-2) circle (3pt);
\filldraw[fill=white] (0,-2) circle (3pt);
\filldraw[fill=white] (1,-3) circle (3pt);
\end{tikzpicture}
}
+
\raisebox{-0.5cm}{
\begin{tikzpicture}[scale=0.3]
\draw (0,1) -- (0,0) -- (-1,-1) -- (-2,-2) -- (-3,-3) -- (-3.5,-3.5);
\draw (0,0) -- (1,-1) -- (1.5,-1.5);
\draw (-1,-1) -- (-0.5,-1.5);
\draw (-3,-3) -- (-2.5,-3.5);
\filldraw (0,0) circle (3pt);
\filldraw (-1,-1) circle (3pt);
\filldraw[fill=white] (-2,-2) circle (3pt);
\filldraw (-3,-3) circle (3pt);
\filldraw[fill=white] (1,-1) circle (3pt);
\end{tikzpicture}
}
+
\raisebox{-0.3cm}{
\begin{tikzpicture}[scale=0.3]
\draw (0,1) -- (0,0) -- (-1,-1) -- (-2,-2) -- (-2.5,-2.5);
\draw (0,0) -- (1,-1) -- (1.5,-1.5);
\draw (-1,-1) -- (0,-2) -- (0.5,-2.5);
\draw (-2,-2) -- (-1.5,-2.5);
\filldraw (0,0) circle (3pt);
\filldraw (-1,-1) circle (3pt);
\filldraw (-2,-2) circle (3pt);
\filldraw[fill=white] (0,-2) circle (3pt);
\filldraw[fill=white] (1,-1) circle (3pt);
\end{tikzpicture}
}.
\end{gather*}
Notice that we can't get the white vertices any deeper in $\mathcal{Q}^2$.
\bigskip

In particular,
\begin{gather} \label{Simply}
\mathcal{Q}^{n+1}\cong X\oo \mathcal{Q}^n.
\end{gather}

Obviously
$$H_*(\oDR^0,\dd^0)\cong \mathcal{Q}^0$$
and we claim that
$$H_*(\oDR^n,\dd^0)\cong \mathcal{Q}^n$$
for $n\geq 1$.
Suppose the claim holds for $n$ and we prove it for $n+1$.
The idea is to use a spectral sequence to get rid of the last sum in the formula
\begin{gather*}
\dd^0(\phi_c^m\oo_1\ul{x}\oo(x_1,\ldots,x_\ar{x})) = \\ \\ \phi_c^{m+1}\oo_1x\oo(x_1,\ldots,x_\ar{x})\ + \\
- \sum_{i=1}^\ar{x}\phi^m\oo_1 x\oo_i \phi \oo (x_1,\ldots,x_\ar{x})\ + \\
+(-1)^{\dg{\ul{x}}}\sum_{i=1}^\ar{x} (-1)^{\sum_{j=1}^{i-1}\dg{x_j}} \phi^m\oo_1 \ul{x}\oo (x_1,\ldots,\dd^0(x_i),\ldots,x_\ar{x}).
\end{gather*}

Consider the spectral sequence $E^{0*}$ on $\oDR^{n+1}$ associated to the filtration
$$0\into\mathfrak{G}_0\into\mathfrak{G}_1\into\cdots\into \oDR^n,$$
where $\mathfrak{G}_k$ is spanned by
$$\phi^m\oo_1 g\oo(x_1,\ldots,x_\ar{x})$$
for all $m\geq 0$, $g\in X\op\ul{X}$, $x_i\in\oDR^n$ and $\sum_{i=1}^\ar{x}\dg{x_i}\leq k$.
Obviously $\dd^0:\mathfrak{G}_k\to \mathfrak{G}_k$.

We will use the comparison theorem for the obvious projection $$\oDR^{n+1}\xrightarrow{\proj} \mathcal{Q}^{n+1}.$$
We consider the zero differential on $\mathcal{Q}^{n+1}$.
It is easily seen that $\proj\dd^0=0$, hence $\proj$ is dg-\coll{C} morphism.
We equip $\mathcal{Q}^{n+1}$ with the trivial filtration $0\into \mathcal{Q}^{n+1}$ and consider the associated spectral sequence $E'^{0*}$.
Again, both filtrations are bounded below and exhaustive.

On the $0^{th}$ page $E^{00}\cong \oDR^{n+1}$, the differential $\dd^{00}$ has the desired form:
\begin{gather*}
\dd^{00}( \phi^m\oo_1\ul{x}\oo(x_1,\ldots,x_\ar{x}) ) = \\
\phi_c^{m+1}\oo_1x\oo(x_1,\ldots,x_\ar{x}) - \sum_{i=1}^\ar{x}\phi^m\oo_1 x\oo_i\phi\oo(x_1,\ldots,x_\ar{x})
\end{gather*}
and $\dd^{00}$ is zero on other elements.
For this differential $\dd^{00}$, it is (at last!) clear how its kernel looks (compare to $\dd^0$), namely $\Ker\dd^{00}=\Fr{}{X\op\Phi}\oo\oDR^n$.
Hence
$$H_*(E^{00},\dd^{00}) \cong X\oo\oDR^n.$$
This is $E^{01}$ and the differential $\dd^{01}$ is equal to the restriction of $\dd^0$ onto $X\oo\oDR^n$.

For $E'^{0*}$ everything is trivial, $E'^{01}\cong \mathcal{Q}^{n+1}$ and $\dd'^{01}=0$.

Then $\proj^1:E^{01}\to E'^{01}$ induced by $\proj$ is quism, because
$$H_*(E^{01},\dd^{01}) \cong X\oo H_*(\oDR^n,\dd^0) \cong X\oo \mathcal{Q}^n \cong \mathcal{Q}^{n+1},$$
where the first isomorphism follows from the usual K\"unneth formula \eqref{UsualKunneth}, the second one follows from the induction hypothesis and the last one was already observed in \eqref{Simply}.

This concludes the proof of the claim $H_*(\oDR^n,\dd^0)\cong \mathcal{Q}^n$.
Finally 
$$H_*(\oDR,\dd^0)\cong H_*(\colim_n\oDR^n)\cong \colim_n H_*(\oDR^n) \cong \colim_n \mathcal{Q}^n \cong \frac{\Fr{}{X\op\Phi}}{\mathfrak{D}}$$
proves Sublemma \ref{sublemma0}.
\end{proof}

Now that the sublemma is proved, we easily go through all the isomorphisms to check \eqref{CheckLater} and \eqref{CheckLater2}.
\end{proof}

\subsection{Augmented cotangent complex} \label{AugmentedCotangentComplex}

Let $$(\oA,\dd_\oA) \xrightarrow{\alg} (\oEnd_A,\dd_{\oEnd_A})$$
be an $\oA$-algebra structure on $A$.
We begin by extracting the operadic cohomology from $\oDR$.
Let $\Fr{}{X\op\Phi\op\ul{X}} \to \oA$ be the dg-\colop{C} morphism which equals $\rR$ on $X$ and vanishes on the other generators.
Hence $\oDR=\Fr{}{X\op\Phi\op\ul{X}}$ is a dg-\colop{C} \emph{over} $\oA$.

For $M$ one of the subsets $X$, $X\op\Phi$, $X\op\ul{X}$ of $\oDR$ define
\begin{gather} \label{DefDerM}
\Der^M_\oA(\oDR,\oEnd_A) := \set{\theta\in\Der_\oA(\oDR,\oEnd_A)}{\forall m\in M\quad \theta(m)=0}.
\end{gather}
We will abbreviate this by $\Der^M$.
Let $\ol{\codd}$ be the differential on $\Der^M$ defined by
$$\ol{\codd}\theta := \theta\ddDR - (-1)^\dg{\theta}\dd_{\oEnd_A}\theta.$$
A check similar to that for \eqref{CotangDiff} verifies this is well defined.
Obviously $$\Der^X = \Der^{X\op\ul{X}} \op \Der^{X\op\Phi}.$$
Recall we assume the dg-$k$-module $A$ is graded by the colours, that is $A=\bigoplus_{c\in C}A_c$.
Hence we have
$$\Der^{X\op\ul{X}} \cong \Hom_{C-\mathrm{coll.}}(\Phi,\oEnd_A) \cong \bigoplus_{c\in C}\Hom_k(A_c,A_c).$$
Importantly, $\Der^{X\op\Phi}$ is closed under $\ol{\codd}$.

\begin{lemma}
$$(\Der^{X\op\Phi},\ol{\codd}) \cong \downar(\Der_\oA(\oR,\oEnd_A),\codd)$$
as dg-$k$-modules.
\end{lemma}

\begin{proof}
Recall $\downarrow\codd=-\codd$.
Define a degree $+1$ map $$\Der^{X\op\Phi}\xrightarrow{f_1} \Der_\oA(\oR,\oEnd_A)$$ by the formula
$$(f_1\theta')(x) := \theta'(\ul{x})\quad\mbox{for }\theta'\in\Der^{X\op\Phi}.$$
Its inverse, $f_2$ of degree $-1$, is defined for $\theta\in\Der_\oA(\oR,\oEnd_A)$ by the formulas
$$(f_2\theta)(a)=0=(f_2\theta)(\phi),\quad (f_2\theta)(\ul{x})=\theta(x).$$
Obviously $f_2f_1=\id$ and $f_1f_2=\id$ and it remains to check $f_1\ol{\codd}=-\codd f_1$.
\begin{align*}
(f_1(\ol{\codd}\theta'))(x) &= (\ol\codd\theta')(\ul{x}) = \theta'(\ddDR\ul{x})-(-1)^{\dg{\theta'}}\dd_{\oEnd_A}(\theta'(\ul{x})), \\
(-\codd(f_1\theta'))(x) &= -(f_1\theta')(\ddR x) + (-1)^{\dg{f_1\theta'}}\dd_{\oEnd_A}((f_1\theta')(x)).
\end{align*}
Now we check $\theta'(\ddDR\ul{x}) = -(f_1\theta')(\ddR x)$.
Let $\ddR x=\sum_i T_i(x_{i1},\ldots,x_{in_i})$.
\begin{align*}
\theta'(\ddDR\ul{x}) &= \theta'(\phi\oo x - \sum_j x \oo_j\phi - s(\ddR x)) = \\
&= -\theta'(s \sum_i T_i(x_{i1},\ldots,x_{in_i})) = \\
&= -\theta'(\sum_i\sum_{j=1}^{n_i}\epsilon_{ij} T_i(x_{i1},\ldots,\ul{x_{ij}},\ldots,x_{in_i})) = \\
&= -\sum_i\sum_j \epsilon_{ij}^{1+\dg{\theta'}} T_i(\rR(x_{i1}),\ldots,\theta'(\ul{x_{ij}}),\ldots,\rR(x_{in_i})), \\
-(f_1\theta')(\ddR x) &= \ldots = -\sum_i\sum_j\epsilon_{ij}^{\dg{f_1\theta'}} T_i( \rR(x_{i1}),\ldots,(f_1\theta')(\ul{x_{ij}}),\ldots,\rR(x_{in_i})  ),
\end{align*}
where we have denoted $\epsilon_{ij}:=(-1)^{\sum_{l=1}^{j-1}\dg{x_{il}}}$.
\end{proof}

\begin{definition}
We call
$$C^*_{aug}(A,A) := ((\Der^X)^{-*}, \ol{\codd})$$
\textbf{augmented operadic cotangent complex} of $A$ and its cohomology
$$H^*_{aug}(A,A) := H^*(C^*_{aug}(A,A), \ol{\codd})$$
\textbf{augmented operadic cohomology} of $A$.
\end{definition}

The interpretation of the augmentation $(\Der^{X\op\ul{X}})^{-*}\xrightarrow{\ol\codd}(\Der^{X\op\Phi})^{-*}\cong C^*(A,A)$ of the usual cotangent complex $C^*(A,A)$ is via infinitesimal automorphisms of the $\oA$-algebra structure on $A$.
This suggests a relation between $H^*_{aug}(A,A)$ and $H^*(A,A)$.
It is best seen in an example:

\begin{example} \label{ExampleAugmentedCohom}
Continuing Example \ref{ExampleDR}, let $A$ be $k$-module with a structure of an associative algebra, that is
$$\oAss\xrightarrow{\alg}\oEnd_A.$$
We have
\begin{align*}
C^0_{aug}(A,A) &= \Der^{X\op\ul{X}} \cong \Hom_k(A,A), \\
C^n_{aug}(A,A) &= (\Der^{X\op\Phi})^{-n} \cong \Hom_{C\mathrm{-coll.}}(\ul{X},\oEnd_A)^{-n} \cong \\
&\cong \Hom_{C\mathrm{-coll.}}(\ul{X}_n,\oEnd_A) \cong \oEnd_A(n+1) = \Hom_k(\otexp{A}{n+1},A)
\end{align*}
and, for $f\in C^n_{aug}(A,A)$,
\begin{gather} \label{DefHochDiff}
\ol{\codd}f = (-1)^{n+1}\mu\oo_2 f + \sum_{k=1}^n (-1)^{n+1-k}f\oo_k\mu + \mu\oo_1 f.
\end{gather}
So the augmented cotangent complex is the Hochschild complex without the term $C^{-1}(A,A)=\Hom_k(k,A)\cong A$,
while the ordinary cotangent complex would be additionally missing $C^0(A,A)$:
$$\underbrace{C^0(A,A) \xrightarrow{\ol{\codd}} \underbrace{C^1(A,A) \xrightarrow{\ol{\codd}} C^2(A,A) \xrightarrow{\ol{\codd}} \cdots}_{C^*(A,A)}}_{C^*_{aug}(A,A)}$$
\end{example}

To generalize the conclusion of the example, recall from \cite{MModels} that $TJ$-grading on a free resolution $\oR=(\Fr{}{X},\dd)\xrightarrow{\rho}(\oA,0)$ is induced by a grading $X=\bigoplus_{i\geq 0}X^i$ on the \coll{C} of generators, denoted by upper indices, $\oR^i$, and satisfying
\begin{enumerate}
\item $\dd$ maps $X^i$ to $\Fr{}{\bigoplus_{j<i}X^j}$,
\item $H_0(\oR^*,\dd) \xrightarrow{H_0(\rho)} \oA$ is an isomorphism of graded \colop{C}s.
\end{enumerate}
If we have a $TJ$-graded resolution $\oR$, we can replace the usual grading by the $TJ$-grading and we let $\Der_\oA(\oR,\oEnd_A)^i$ be the $k$-module of derivations $\oR\to \oEnd_A$ vanishing on all $X^j$'s except for $j=i$ and let ${}^{TJ}C^{*}(A,A):=\upar(\Der_\oA(\oR,\oEnd_A)^{*},\codd)$ and
${}^{TJ}H^{*}(A,A):=H^*({}^{TJ}C^{*}(A,A))$.
In case $\oA$ is concentrated in degree $0$, the usual grading is $TJ$ and we get the same result as in \eqref{CotangDef}, i.e. $C^*(A,A)={}^{TJ}C^{*}(A,A)$ and we can forget about the superscripts $TJ$ everywhere.

For a $TJ$-graded $\oR$, we can also equip $C^*_{aug}(A,A)$ with similar $TJ$-grading ${}^{TJ}C^*_{aug}(A,A)$ as above.
On this matter we just remark that $\phi$ is placed in $TJ$-degree $0$ and leave the details for the interested reader.
Finally, the following is obvious:

\begin{theorem}
\begin{enumerate}
\item ${}^{TJ}C^n_{aug}(A,A) = 0$ for $n\leq -1$,
\item ${}^{TJ}C^0_{aug}(A,A) = \Der^{\oA\op\ul{X}} \cong \Hom_k(A,A)$,
\item ${}^{TJ}H^1_{aug}(A,A) \cong k$-module of formal infinitesimal deformations of the $\oA$-algebra structure on $A$ \emph{modulo infinitesimal automorphisms},
\item ${}^{TJ}H^n_{aug}(A,A) \cong {}^{TJ}H^n(A,A)$ for $n\geq 2$.
\end{enumerate}
\end{theorem}

Notice that the unaugmented operadic cohomology ${}^{TJ}H^1(A,A)$ is the $k$-module of formal infinitesimal deformations of the $\oA$-algebra structure on $A$, but the infinitesimal automorphisms are not considered.

Hence the distinction between $H^*(A,A)$ and $H^*_{aug}(A,A)$ is inessential and we will usually not distinguish these two.

\subsection{Intermediate resolution of $\oDA$} \label{SectionIntermediate}

Now we construct an intermediate step in the resolution of Theorem \ref{ResolutionOfDerivation}:
\begin{diagram}
\oDR & & \rTo^\sim_{\rDR} & & \oDA \\
 & \rdTo^\sim_\interm & & \ruTo^\sim_\rolDR \\
 & & \ol{\oDR} & &
\end{diagram}
Intuitively, $\iota$ should "unresolve" the part of $\oDR$ corresponding to the $\oA$-algebra operations and do nothing in the part corresponding to the derivation $\phi$.
Let
$$\olDR := \left( \oA \fpr \Fr{}{\Phi\op\ul{X}} , \ddolDR \right).$$
We first define $\interm$ to be the composite
$$\oDR = \Fr{}{X\op\Phi\op\ul{X}} \cong \Fr{}{X} \fpr \Fr{}{\Phi\op\ul{X}} \xrightarrow{\rR \fpr \id} \oA \fpr \Fr{}{\Phi\op\ul{X}} = \olDR$$
then $\ddolDR$ is the derivation defined by 
\begin{align}
&\ddolDR (a) := \ddDA(a) = \ddA(a), \nonumber \\
&\ddolDR (\phi_c) := 0, \label{DefOfddolDR} \\
&\ddolDR (\ul{x}) :=  \interm(\ddDR\ul{x}). \nonumber
\end{align}
Now we check $\interm \ddDR = \ddolDR \interm$ and this will immediately imply $\ddolDR^2=0$:
\begin{gather*}
\interm \ddDR (x) = \interm \ddR (x) = \rR \ddR (x) = \ddA \rR (x), \\
\ddolDR \interm (x) = \ddolDR \rR(x) = \ddDA \rR(x) = \ddA \rR (x)
\end{gather*}
and similar claim for $\ul{x}$ is an immediate consequence of definitions.

Finally, let $\rolDR$ be the \colop{C} morphism defined by
\begin{align}
&\rolDR(a) := a, \nonumber \\ 
&\rolDR(\phi_c) := \phi_c, \label{DefOfrolDR} \\ 
&\rolDR(\ul{x}) := 0. \nonumber
\end{align}

\begin{lemma}
$\rolDR$ is dg-\colop{C} morphism.
\end{lemma}

\begin{proof}
We only have to check $\rolDR\ddolDR(\ul{x}) = 0$:
\begin{align*}
\rolDR\ddolDR(\ul{x}) &= \rolDR\interm\left( \phi_c \oo_1 x - \sum_{i=0}^n x\oo_i \phi_{c_i} - s(\ddR x) \right) = \\
&= \phi_c \oo_1 \rR(x) - \sum_{i=0}^n \rR(x)\oo_i \phi_{c_i}.
\end{align*}
The third term in the bracket vanishes since $\interm s \ddR (x)$ is a sum of compositions each of which contains a generator from $\ul{X}$ and $\rolDR$ vanishes on $\ul{X}$.
The above expression vanishes because if $\rR(x)\neq 0$, then it is precisely the relator \eqref{relator}.
\end{proof}

\begin{lemma}
$\interm$ is a quism.
\end{lemma}

\begin{proof}
We notice that $\olDR$ is close to be the $1^{st}$ term of a spectral sequence computing homology of $\oDR$.
Now we make this idea precise.

Consider a new grading $\gr$ on $\oDR$:
$$\gr(x):=0=:\gr(\phi),\quad \gr(\ul{x}):=\dg{\ul{x}}$$
and its associated filtration
\begin{gather*}
\mathfrak{F}_p:=\bigoplus_{i=0}^p\set{z\in\oDR}{\gr(z)=i}, \\
0\into\mathfrak{F}_0\into\mathfrak{F}_1\into\cdots,\quad \ddDR\mathfrak{F}_p\subset\mathfrak{F}_p
\end{gather*}
and its associated spectral sequence $(E^*,\dd^*)$.
There is an analogous spectral sequence $(E'^*,\dd'^*)$ on $\olDR$ given by the grading
$$\gr'(a):=0=:\gr'(\phi),\quad \gr'(\ul{x}):=\dg{\ul{x}}.$$
Since both filtrations are bounded below and exhaustive, we can use the comparison theorem.

We have $E^0\cong \Fr{}{X\op\Phi\op\ul{X}}$.
Recalling the formulas \eqref{DiffOnoDR}, we immediately see that $\dd^0$ on $E^0$ is the derivation differential given by
$$\dd^0x=\ddR x,\quad \dd^0\phi=0=\dd^0\ul{x}.$$
Hence $E^1\cong H_*(E^0,\dd^0)\cong H_*(\Fr{}{X})\fpr H_*(\Fr{}{\Phi\op \ul{X}}) \cong\oA\fpr\Fr{}{\Phi\op\ul{X}}$ by the K\"unneth formula for a free product of dg-\colop{C}s, see Lemma \ref{LemmaKunneth}.
Similarly $E'^1\cong \oA\fpr\Fr{}{\Phi\op\ul{X}}$.

Understanding the differentials $\dd^1$ and $\dd'^1$ on the $1^{st}$ pages as well as the induced dg-\coll{C} morphism $\interm^1$ is easy (though notationally difficult - observe $\dd^1\ul{x}$ is \emph{not} $\interm(\ddDR\ul{x})$ in general!) and we immediately see that $\interm^1$ is an isomorphism of dg-\coll{C}s.
\end{proof}

\begin{corollary} \label{FirstCorollary}
$\rolDR$ is a resolution of $\oDA$.
\end{corollary}

\begin{example} \label{ExampleolDR}
Let's continue Example \ref{ExampleDR} and make $\olDR$ explicit:
$$\olDR:=(\oAss\fpr\Fr{}{\Phi\op\ul{X}},\ddolDR),$$
\begin{align*}
&\ddolDR(a) := 0 =: \ddolDR(\phi), \\
&\ddolDR(\ul{x}^2) := \phi\oo\mu - \mu\oo_1\phi - \mu\oo_2\phi, \\
&\ddolDR(\ul{x}^n) := -(-1)^n \mu\oo_2 \ul{x}^{n-1} - \sum_{k=1}^{n-1} (-1)^{n-k} \ul{x}^{n-1}\oo_k\mu - \mu\oo_1\ul{x}^{n-1}
\end{align*}
for $n\geq 3$.
The last formula is reminiscent to the one for the Hochschild differential.
We will make this point precise in Section \ref{FromOperadsToModules}.
\end{example}

\subsection{From operads to operadic modules} \label{FromOperadsToModules}

Associated to the operad $\oDA = (\oA \fpr \Fr{}{\Phi}) / \mathfrak{D}$ is the dg-$(\oA,\ddA)$-module
$$\oMDA := \left( \frac{\FrM{\oA}{\Phi}}{\mathfrak{D} \cap \FrM{\oA}{\Phi}} , \ddMDA \right),$$
where $\mathfrak{D} \cap \FrM{\oA}{\Phi}$ is the sub-$\oA$-module of $\FrM{\oA}{\Phi}$ generated by the relators \eqref{relator} and $\ddMDA$ is a dg-$\oA$-module morphism given by $\ddMDA\phi=0$.

Associated to $\olDR = ( \oA\fpr\Fr{}{\Phi\op\ul{X}} , \ddolDR )$ is the dg-$\oA$-module
$$\olMDR := \left( \FrM{\oA}{\Phi\op\ul{X}} , \ddolMDR \right),$$
where $\ddolMDR$ is a $\oA$-module morphism given by the same formulas \eqref{DefOfddolDR} as $\ddolDR$.
We emphasize that this makes sense because $\FrM{\oA}{\Phi\op\ul{X}} \subset \oA\fpr\Fr{}{\Phi\op\ul{X}}$, the dg-$\oA$-module structure is induced by the operadic composition and $\ddolDR$ maps $\FrM{\oA}{\Phi\op\ul{X}}$ into itself!

Associated to the dg-\colop{C} morphism $\rolDR : \olDR \to \oDA$ is the dg-$\oA$-module morphism
$$\rolMDR : \olMDR \to \oMDA$$
again defined by the formulas \eqref{DefOfrolDR} as $\rolMDR$.

\begin{lemma} \label{LemmaAboutrolDR}
$\rolMDR$ is a quism.
\end{lemma}

\begin{proof}
Let $G_p$ be the sub-\coll{C} of $\olDR = \oA\fpr\Fr{}{\Phi\op\ul{X}}$ spanned by the compositions containing precisely $p$ generators from $\Phi\op\ul{X}$, i.e. $G_0=\oA$, $G_1=\FrM{\oA}{\Phi\op\ul{X}}$ and 
$$\olDR = \bigoplus_{p\geq 0} G_p.$$
We have analogous grading $\oA\fpr\Fr{}{\Phi} = \bigoplus_{p\geq 0} G'_p$.
Let $\proj:\oA\fpr\Fr{}{\Phi} \to \oDA$ be the natural projection.
Since relators \eqref{relator} are homogeneous with respect to this grading, $G''_p:=\proj G'_p$ defines a grading
$$\oDA = \bigoplus_{p\geq 0}G''_p.$$

Observe $G''_0=\oA$ and $G''_1=\oMDA$.
By definitions, $\rolDR G_p \subset G''_p$, hence $\rolDR$ decomposes as a sum of $\rolDR^p: G_p \to G''_p$.
The above direct sums are in fact direct sums of sub-dg-\coll{C}s, $\rolDR$ is a quism by Corollary \ref{FirstCorollary}, hence all the $\rolDR^p$'s are quisms, especially $\rolDR^1=\rolMDR$.
\end{proof}

Now we formalize the statement : 
\emph{$\olMDR$ contains all the information needed to construct the operadic cohomology for $\oA$-algebras.}

First observe that $\oEnd_A$ is naturally a dg-$\oA$-module.
Let $\ol{\codd}$ be the differential on $\Hom_{\dgmod{\oA}}(\olMDR,\oEnd_A)$ defined by the formula 
\begin{gather*}
\ol{\codd}\theta:=\theta\dd_\olMDR - (-1)^\dg{\theta}\dd_{\oEnd_A}\theta
\end{gather*}
similar to \eqref{CotangDiff}.

\begin{lemma}
$$(C^*_{aug}(A,A),\ol{\codd}) \cong \left( \Hom_{\dgmod{\oA}}(\olMDR,\oEnd_A),\ol\codd \right)$$
as dg-$k$-modules.
\end{lemma}

\begin{proof}
On the level of $k$-modules, we have
\begin{gather*}
C^*_{aug}(A,A) = \Der^X = \set{\theta\in\Der_\oA(\fpr\Fr{}{X\op\Phi\op\ul{X}},\oEnd_A)}{\forall x\in X\quad \theta(x)=0} \cong \\
\cong \Hom_{\mathrm{dg-}C\mathrm{-coll.}}(\Phi\op\ul{X},\oEnd_A)
\end{gather*}
by the defining property of derivations and
$$\Hom_{\dgmod{\oA}}(\FrM{\oA}{\Phi\op\ul{X}},\oEnd_A) \cong \Hom_{\mathrm{dg-}C\mathrm{-coll.}}(\Phi\op\ul{X},\oEnd_A)$$
by the freeness of $\FrM{\oA}{\Phi\op\ul{X}}$.

The differentials are clearly preserved under the above isomorphism.
\end{proof}

The nice thing is that we have now all the information needed to construct the cohomology for $\oA$-algebras encoded in terms of the \emph{abelian} category of dg-$\oA$-modules. Hence

\begin{theorem} \label{CohomAsExt}
\begin{gather*}
H^*_{aug}(A,A) \cong H_*( \Hom_{\dgmod{\oA}}^{-*}(\olMDR,\oEnd_A),\ol\codd ) \cong \\
\cong \Ext^{-*}_{\dgmod{\oA}}(\oMDA,\oEnd_A)
\end{gather*}
\end{theorem}

In particular, since the homotopy theory in abelian categories is well known and simple, it is immediate that $\Ext_{\dgmod{\oA}}$ and hence $H^*_{aug}(A,A)$ doesn't depend on the choice of a projective resolution of $\oMDA$ and consequently doesn't depend on the choice of the free resolution $\oR\xrightarrow[\rR]{\sim}\oA$ in \eqref{resolutionR}.

The main advantage of the above expression is that in order to construct cohomology for $\oA$-algebras, we don't need to find a free (or cofibrant) resolution $\oR\xrightarrow{\sim}\oA$ in the category of dg-\colop{C}s, but it suffices to find a projective resolution of $\oMDA$ in the category of dg-$\oA$-modules, which is certainly easier.

\begin{example} \label{ExampleAssModule}
Let's continue Example \ref{ExampleolDR}:
$$\oMDA:=\frac{\FrM{\oAss}{\Phi}}{(\phi\oo\mu-\mu\oo_1\phi-\mu\oo_2\phi)}$$
and we have the following explicit description of $\olMDR$:
$$\olMDR := (\FrM{\oAss}{\kspan{\phi^1,\phi^2,\phi^3,\ldots}},\ddolMDR),$$
where $\phi^1:=\phi$ and $\phi^n:=\ul{x}^n$, for $n\geq 2$, is of degree $n-1$ and the differential is given by
\begin{align*}
&\ddolMDR(\phi^1) := 0, \\
&\ddolMDR(\phi^n) := -(-1)^n \mu\oo_2 \phi^{n-1} - \sum_{k=1}^{n-1} (-1)^{n-k} \phi^{n-1}\oo_k\mu - \mu\oo_1\phi^{n-1}.
\end{align*}
Lemma \ref{LemmaAboutrolDR} states that $\olMDR \xrightarrow{\rolMDR} \oMDA$ is a free resolution in the category of $\oAss$-modules.

Notice the similarity to the Hochschild complex.
This suggests that if we know a complex computing a cohomology for $\oA$-algebras (in this case Hochschild complex) we can read off a candidate for the free resolution of $\oMDA$ (in this case we already know a resolution $\oMDA$, namely $\olMDR$, but this was constructed from the \emph{operadic} resolution $\oR$, which is not generally available).
If we can prove that this candidate is indeed a resolution, we get that the cohomology in question is isomorphic to the augmented operadic cohomology.
\end{example}

We demonstrate this by constructing a cohomology for diagrams of associative algebras and proving that the (augmented) cotangent complex coincides with that defined by Gerstenhaber and Schack \cite{GS}.

On the other hand, in the process of constructing $\olMDR$ we have discarded much information present in $\oR$. Namely $\oR$ can be used to define an $L_\infty$ structure on $C^*(A,A)$ governing formal deformations of $A$ (see \cite{IB}), which is no longer possible using $\olMDR$ (or any other resolution of $\oMDA$) only.

\section{Gerstenhaber-Schack diagram cohomology is operadic cohomology} \label{SectionGS}

\subsection{Operad for diagrams}

Let $\cC$ be a small category.
For a morphism $f$ of $\cC$, let $\inp{f}$ be its source (Input) and $\out{f}$ its target (Output).
Consider the following nerve construction on $\cC$:
$$\Sigma^n := \set{ \left(\xleftarrow{f_n}\cdots\xleftarrow{f_1}\right)\in\Hom_\cC^{\times n} }{\out{f_i}=\inp{f_{i+1}}\mbox{ for }1\leq i\leq n-1}$$
for $n\geq 1$.
For $\sigma=\left(\xleftarrow{f_n}\cdots\xleftarrow{f_1}\right)$, let $\dg{\sigma}:=n$, let $\inp{\sigma}:=\inp{f_1}$ and $\out{\sigma}:=\out{f_n}$.
The face maps $\Sigma^n\to\Sigma^{n+1}$ are given by $\sigma\mapsto \sigma_i$, where
\begin{align*}
\sigma_0 &:= \left(\xleftarrow{f_n}\cdots\xleftarrow{f_2}\right), \\
\sigma_i &:= \left(\xleftarrow{f_n}\cdots\xleftarrow{f_{i+1}f_i}\cdots\xleftarrow{f_1}\right)\mbox{ for }1\leq i\leq n-1, \\
\sigma_n &:= \left(\xleftarrow{f_{n-1}}\cdots\xleftarrow{f_1}\right).
\end{align*}
Denote $\Sigma^0$ the set of objects of $\cC$ and for $\sigma\in\Sigma^0$, let $\inp{\sigma}=\out{\sigma}=\sigma$.
Finally let
$$\Sigma:=\bigcup_{n=0}^\infty \Sigma^n$$
and denote $\Sigma^{\geq 1}:=\Sigma-\Sigma^0$.

Let $\oC$ be the operadic version of $\cC$, that is
$$\oC:=\kspan{\Sigma^1}.$$
This can be seen as a \colop{\Sigma^0}, where each $f\in\Sigma^1$ is an element of $\oC\binom{\out{f}}{\inp{f}}$ and the operadic composition is induced by the categorical composition.

A ($\cC$-shaped) \textbf{diagram} (of associative algebras) is a functor $$D:\cC\to\oAss\mbox{-algebras}.$$
Now we describe a \colop{\Sigma^0} $\oA$ such that $\oA$-algebras are precisely $\cC$-shaped diagrams:
$$\oA:=\frac{( \bigfpr_{c\in\Sigma^0} \oAss_c)\fpr\oC}{\mathcal{I}},$$
where $\oAss_c$ is a copy of $\oAss$ concentrated in colour $c$, its generating element is $\mu_c\in\oAss_c\binom{c}{c,c}$ and $\mathcal{I}$ is the ideal generated by
$$f\oo\mu_{\inp{\sigma}}-\mu_{\out{\sigma}}\oo(f,f)\quad\mbox{for all }f\in\Sigma^1.$$
It should be clear now that the functor $D$ is essentially the same thing as \colop{\Sigma^0} morphism $\oA\to\oEnd_A$, where $A=\bigoplus_{c\in\Sigma^0}D(c)$.

The associated module of Section \ref{FromOperadsToModules} is 
$$\oMDA:=\frac{\FrM{\oA}{\bigoplus_{c\in\Sigma^0}\Phi_c}} {\mathfrak{D}\cap\FrM{\oA} {\bigoplus_{c\in\Sigma^0}\Phi_c}}$$
where $\Phi_c=\kspan{\phi_c}$, $\phi_c$ being an element of colour $\binom{c}{c}$ and of degree $0$, and 
the submodule in the denominator is generated by
\begin{gather*}
\phi_c\oo\mu_c - \mu_c\oo_1\phi_c - \mu_c\oo_2\phi_c, \\
\phi_{\out{f}}\oo f - f\oo\phi_{\inp{f}}
\end{gather*}
for all $c\in\Sigma^0$ and all $f\in\oC$ (equivalently $f\in\Sigma^1$).
We seek a free resolution $(\oMR,\dd)\xrightarrow{\sim}(\oMDA,0)$ to use Theorem \ref{CohomAsExt}.
Before constructing $\oMR$, let's recall the Gerstenhaber-Schack diagram cohomology.
As we have seen in Example \ref{ExampleAssModule}, this gives us a candidate for $\oMR$.

\subsection{Gerstenhaber-Schack diagram cohomology}

We adapt the notation from the original source \cite{GS}.
Originally, the diagram $D$ was restricted to be a poset, but this is unnecessary.
Also, instead of associative algebras, one may consider any other type of algebras for which a convenient cohomology is known (e.g. Lie algebras, \cite{GGS}).
In this paper we stick to associative algebras, but we believe that other types can be handled in a similar way.

For $\sigma\in\Sigma^0$, denote $\ul{\sigma}:=\id_\sigma$.
For $\sigma=( \xleftarrow{f_p}\cdots\xleftarrow{f_1} ) \in\Sigma^p$, denote
$$\ul{\sigma}:=f_p\cdots f_1 : D(\inp{\sigma})\to D(\out{\sigma})$$ the composition along $\sigma$.
This algebra morphisms makes $D(\out{\sigma})$ a $D(\inp{\sigma})$-bimodule.
For $p,q\geq 0$ let
$$C^{p,q}_{GS}(D,D):=\prod_{\sigma\in\Sigma^p} C^q_{\mathrm{Hoch}}(D(\inp{\sigma}),D(\out{\sigma})),$$
where $C^q_{\mathrm{Hoch}}(D(\inp{\sigma}),D(\out{\sigma})) = \Hom_{k}(\otexp{D(\inp{\sigma})}{q+1},D(\out{\sigma}))$ are the usual Hochschild cochains.
We usually abbreviate $C^{p,q}_{GS}:=C^{p,q}_{GS}(D,D)$.
There are vertical and horizontal differentials $\codd_V : C^{p,q}_{GS} \to C^{p,q+1}_{GS}$ and $\codd_H:C^{p,q}_{GS} \to C^{p+1,q}_{GS}$.
To write them down, let $\sigma=( \xleftarrow{f_{p+1}}\cdots\xleftarrow{f_1} )\in\Sigma^{p+1}$ and for $\tau\in\Sigma^p$ let $\proj_\tau : \prod_{\lambda\in\Sigma^p}C^q_{\mathrm{Hoch}}(D(\inp{\lambda}),D(\out{\lambda})) \to C^q_{\mathrm{Hoch}}(D(\inp{\tau}),D(\out{\tau}))$ be the projection onto the $\tau$ component of $C^{p,q}_{GS}$.
Let $\coddH$ be the usual Hochschild differential, see \eqref{DefHochDiff}.
Finally, for $\theta\in C^{p,q}_{GS}$, let
\begin{align}
&\codd_V := (-1)^p\prod_{\sigma\in\Sigma^p} \coddH, \label{coddVertical} \\
&\proj_\sigma(\codd_H\theta) := (-1)^{p+1}  (\proj_{\sigma_0}\theta)\oo(\underbrace{f_1,\ldots,f_1}_{\scriptstyle{q\mathrm{-times}}}) 
+ \sum_{i=1}^p (-1)^{p+1-i} \proj_{\sigma_i}\theta \ + \nonumber \\
&\phantom{\proj_\sigma(\codd_H\theta) :=} + f_{p+1}\oo(\proj_{\sigma_{p+1}}\theta). \label{coddHorizontal}
\end{align}

It is easy to see that $(C^{*,*}_{GS},\codd_V,\codd_H )$ is a bicomplex.
The Gerstenhaber-Schack cohomology is defined to be the cohomology of the totalization of this bicomplex,
$$H^*_{GS}(D,D) := H^*(\bigoplus_{p+q=*}C^{p,q}_{GS}(D,D),\codd_V+\codd_H).$$

Notice that we have restricted ourselves to the Hochschild complex without $C^{-1}_{\mathrm{Hoch}}$ as in Example \ref{ExampleAugmentedCohom}.
Also we consider only the cohomology of $D$ with \emph{coefficients in itself} as this is the case of interest in the formal deformation theory.
The general coefficients can be handled using trivial (operadic) extensions.

\subsection{Resolution of $\oMDA$} \label{SectionResolution}

The preceding section gives us the following candidate for $\oMR$:
$$\oMR:=\left( \FrM{\oA}{\bigoplus_{\sigma\in\Sigma}\Phi_\sigma\op X_\sigma},\dd \right),$$
where $\Phi_\sigma:=\kspan{\phi_\sigma}$ with $\phi_\sigma$ of colours $\binom{\out{\sigma}}{\inp{\sigma}}$ and of degree $\dg{\phi_\sigma}:=\dg{\sigma}$ and
$X_\sigma:=\upar^{\dg{\sigma}+1}X$ is placed in output colour $\out{\sigma}$ and input colours $\inp{\sigma}$, $X$ being the collection of generators of the minimal resolution of $\oAss$ as in Example \ref{ExampleDR}.
The element of $X_\sigma$ corresponding to $x\in X$ will be denoted by $x_\sigma$, hence $\dg{x^i_\sigma}=i-1+\dg{\sigma}$.
To define the differential $\dd$ in an economic way, denote $x^1_\sigma:=\phi_\sigma$ for $\sigma\in\Sigma^0$ and let
$$\pre{\phi_\sigma}:=0,\quad \pre{x^i_\sigma}:=x^{i-1}_\sigma\mbox{ for }i\geq 2,$$
and extend linearly to the generators of $\oMR$.
Further, let's accept the convention that for $\sigma\in\Sigma^0$, the symbol $x_{\sigma_0}$ stands for zero.
Then 
\begin{align*}
&\dd(x_\sigma) := (-1)^{\dg{\sigma}}\left( (-1)^{\dg{x}}\mu_{\out{\sigma}}\oo(\ul{\sigma},\pre{x_\sigma}) \ +\phantom{\sum_{i=1}^{\dg{x}+1}}\right. \\ 
&\phantom{\dd(x_\sigma) :=}\left. + \sum_{i=1}^{\ar{x}-1}(-1)^{\dg{x}-i}\pre{x_\sigma}\oo_i\mu_{\inp{\sigma}} + \mu_{\out{\sigma}}\oo(\pre{x_\sigma},\ul{\sigma}) \right) + \\
&\phantom{\dd(x_\sigma) :=}+ (-1)^{\dg{\sigma}}x_{\sigma_0}\oo(\underbrace{f_1,\ldots,f_1}_{(\ar{x})\textrm{-times}}) + \sum_{i=1}^{\dg{\sigma}-1}(-1)^{\dg{\sigma}-i}x_{\sigma_i} + f_{\dg{\sigma}}\oo x_{\sigma_{\dg{\sigma}}}
\end{align*}
for any $\sigma\in\Sigma$.
Observe that the first part of the above formula corresponds to the vertical differential \eqref{coddVertical} and the second part corresponds to the the horizontal differential \eqref{coddHorizontal} in $C^*_{GS}(D,D)$.
Then it is easily seen that $(\Hom_{\dgmod{\oA}}(\oMR,\oEnd_{\bigoplus_{c\in\Sigma^0}D(c)}),\codd)$ with $\codd(-):=-\oo\dd$ is, as a dg-$k$-module, isomorphic to the Gerstenhaber-Schack complex.
Once we prove that $\oMR$ is a resolution of $\oMDA$, we will have, by Theorem \ref{CohomAsExt},

\begin{theorem}
Gerstenhaber-Schack diagram cohomology $H^*_{GS}(D,D)$ is isomorphic to the augmented operadic cohomology $H_{aug}^*(D,D)$.
\end{theorem}

To prove that $\oMR$ is a resolution of $\oMDA$ we introduce the dg-$\oA$-module morphism $\rho:(\oMR,\dd)\to(\oMDA,0)$ given by the formulas
$$\begin{array}{l}
\rho(\phi_\sigma) := \left\{ \begin{array}{lcl} \phi_\sigma & \ldots & \dg{\sigma}=0 \\
																								0						& \ldots & \dg{\sigma}\geq 1,
												 		 \end{array} \right. \\
\rho(x_\sigma) := 0. \\
\end{array}$$
Indeed, it is easy to check that $\rho\dd=0$.
It remains to prove

\begin{lemma} \label{RhoIsQuism}
$\rho$ is a quism.
\end{lemma}

\begin{proof}
is basically a reduction to the following two cases:
\begin{enumerate}
\item $\cC$ is a single object with no morphism except for the identity (Lemma \ref{LemmaA}),
\item $\cC$ is arbitrary, but each $D(c)$, $c\in\Sigma^0$, is the trivial algebra $k$ with zero multiplication (Lemma \ref{LemmaB}).
\end{enumerate}
We first give a general overview of the proof and postpone technicalities to subsequent lemmas.

Consider a new grading on $\oMR$ given by
$$\gr(x_\sigma):=\dg{\sigma}=:\gr(\phi_\sigma)$$
and the usual requirement that the composition is of degree $0$.
Then we have the associated filtration 
\begin{gather*}
\mathfrak{F}_n:=\bigoplus_{i=0}^n\set{x\in\oMR}{\gr(x)=i}, \\
0\into \mathfrak{F}_0 \into \mathfrak{F}_1 \into\cdots,\quad \dd\mathfrak{F}_i\subset \mathfrak{F}_i
\end{gather*}
and the spectral sequence $(E^*,\dd^*)$ which is convergent as the filtration is bounded below and exhaustive.

Obviously $E^0\cong\oMR$ and $\dd^0$ is the derivation differential given by
\begin{align}
&\dd^0(x_\sigma) = (-1)^{\dg{\sigma}}\left( (-1)^{\dg{x}}\mu_{\out{\sigma}}\oo(\ul{\sigma},\pre{x_\sigma}) \ +\phantom{\sum_{i=1}^{\dg{x}+1}}\right. \label{Defdd0} \\ 
&\phantom{\dd^0(x_\sigma) = }\left. + \sum_{i=1}^{\ar{x}-1}(-1)^{\dg{x}-i}\pre{x_\sigma}\oo_i\mu_{\inp{\sigma}} + \mu_{\out{\sigma}}\oo(\pre{x_\sigma},\ul{\sigma}) \right), \nonumber \\
&\dd^0(\phi_\sigma) = 0. \nonumber
\end{align}
Now $H_*(E^0,\dd^0)\cong\bigoplus_{\sigma\in\Sigma}H_*(\FrM{\oA}{\Phi_\sigma\op X_\sigma},\dd^0)$ and we use

\begin{lemma} \label{LemmaA}
$$H_*(\FrM{\oA}{\Phi_\sigma\op X_\sigma},\dd^0) \cong \frac{\FrM{\oA}{\Phi_\sigma}}{\mathfrak{D}_\sigma},$$
where $\mathfrak{D}_\sigma$ is the submodule generated by
\begin{gather} \label{FormulaOne}
\mu_{\out{\sigma}}\oo(\ul{\sigma},\phi_\sigma) + \mu_{\out{\sigma}}\oo(\phi_\sigma,\ul{\sigma}) - \phi_\sigma\oo\mu_{\inp{\sigma}}.
\end{gather}
\end{lemma}

This lemma implies
\begin{align}
&E^1 \cong \frac{\FrM{\oA}{\bigoplus_{\sigma\in\Sigma}\Phi_\sigma}}{\bigoplus_{\sigma\in\Sigma}\mathfrak{D}_\sigma}, \nonumber \\
&\dd^1(\phi_\sigma) = (-1)^{\dg{\sigma}}\phi_{\sigma_0}\oo f_1 + \sum_{i=1}^{\dg{\sigma}-1}(-1)^{\dg{\sigma}-i}\phi_{\sigma_i} + f_{\dg{\sigma}}\oo\phi_{\sigma_{\dg{\sigma}}}. \label{FormulaBarCobar}
\end{align}
Denote $\Phitot:=\bigoplus_{\sigma\in\Sigma}\Phi_\sigma$ and write the nominator in the form $$\FrM{\oA}{\Phitot} \cong \oA\oo'(I,\Phitot\oo\oA).$$
Here we used the infinitesimal composition product \eqref{ICP}.
Because of the relations $f\oo\mu_{\inp{f}}-\mu_{\out{f}}\oo(f,f)$ in $\oA$ for all $f\in\Sigma^1$, we have
\begin{gather} \label{DerivationRelationsUsed}
\oA \cong (\bigoplus_{c\in\Sigma^0}\oAss_c)\oo\oC
\end{gather}
and hence $$\FrM{\oA}{\Phitot} \cong \oA\oo'(I,\Phitot\oo(\bigoplus_{c\in\Sigma^0}\oAss_c)\oo\oC).$$
But we are interested in the quotient $E^1$ of $\FrM{\oA}{\Phitot}$ and the corresponding relations \eqref{FormulaOne} give us
$$E^1 \cong \oA\oo'(I,\Phitot\oo\oC).$$
Now we use \eqref{DerivationRelationsUsed} again to obtain
$$E^1 \cong (\bigoplus_{c\in\Sigma^0}\oAss_c) \oo' (\oC,\oC\oo\Phitot\oo\oC).$$
Notice that $\oC\oo\Phitot\oo\oC \cong \FrM{\oC}{\Phitot}$.
Since $\dd^1$ is nontrivial only on $\Phitot$ and $\FrM{\oC}{\Phitot}$ is closed under $\dd^1$, to understand $H_*(E^1,\dd^1)$ using the usual K\"unneth formula, we only have to compute

\begin{lemma} \label{LemmaB}
$$H_*(\FrM{\oC}{\Phitot},\dd^1) \cong \frac{\FrM{\oC}{\bigoplus_{c\in\Sigma^0}\Phi_c}}{\mathfrak{D}'},$$
where $\mathfrak{D}'$ is the submodule generated by
\begin{gather} \label{FormulaTwo}
f\oo\phi_{\inp{f}}-\phi_{\out{f}}\oo f
\end{gather}
for all $f\in\oC$.
\end{lemma}

Then, tracking back all the above isomorphisms, we get
$$E^2 \cong H_*(E^1,\dd^1) \cong \frac{\FrM{\oA}{\bigoplus_{c\in\Sigma^0}\Phi_c}}{\mathfrak{D}''},$$
where $\mathfrak{D}''$ is the submodule generated by relators \eqref{FormulaOne} for $\sigma\in\Sigma^0$ and all \eqref{FormulaTwo}'s.
Hence $E^2 \cong \oMDA$ and this is concentrated in degree $0$, the spectral sequence collapses and this concludes the proof of Lemma \ref{RhoIsQuism}.
\end{proof}

\begin{proof}[of Lemma \ref{LemmaA}]
We have already seen in Example \ref{ExampleAssModule} that for $c\in\Sigma^0$ the restriction of $\rho$,
\begin{gather} \label{easyFUJ}
\left( \FrM{\oAss_c}{\Phi_c\op X_c},\dd^0 \right) \xrightarrow{\rho} \frac{\FrM{\oAss_c}{\Phi_c}}{\mathfrak{J}_c},
\end{gather}
is a quism, where $\mathfrak{J}_c$ is the submodule generated by \eqref{FormulaOne} for $\sigma=c$.
We will reduce our problem to this case.
Let $$\oM_\sigma := \oAss_{\out{\sigma}}\oo'(\kspan{\ul{\sigma}},(\Phi_\sigma\op X_\sigma)\oo\oAss_{\inp{\sigma}}).$$
This is in fact a sub-\coll{\Sigma^0} of $\FrM{\oA}{\Phi_\sigma\op X_\sigma}$.
An easy computation shows that it is closed under $\dd^0$.
It will play a role similar to $\FrM{\oAss_c}{\Phi_c\op X_c}$ above:

\begin{sublemma} \label{SublemmaFUJ}
There is an isomorphism
$$(\oM_\sigma,\dd^0) \cong \upar^{\dg{\sigma}}\left( \FrM{\oAss_c}{\Phi_c\op X_c} , \dd^0 \right)$$
of dg\coll{}s (we \emph{ignore the colours!}).
This induces an isomorphism
$$H_*(\oM_\sigma,\dd^0) \cong \frac{\oAss_{\out{\sigma}}\oo'(\kspan{\ul{\sigma}},\Phi_\sigma\oo\oAss_{\inp{\sigma}})} {\mathfrak{J}_\sigma}$$
of \coll{\Sigma^0}s (compare to the right-hand side of \eqref{easyFUJ}), where the quotient by $\mathfrak{J}_\sigma$ expresses the fact that $\phi_\sigma$ behaves like a derivation with respect to $\mu_{\out{\sigma}}$ and $\mu_{\inp{\sigma}}$ in $\oM_\sigma$.
The relators in $\mathfrak{J}_\sigma$ are analogous to \eqref{FormulaOne}, namely $\mathfrak{J}_\sigma$ is sub-\coll{\Sigma^0} of $\oM_\sigma$ consisting of elements
\begin{gather*}
(a_{\out{\sigma}}\oo_i(\mu_{\out{\sigma}}\oo_1 a^1_{\out{\sigma}})) \oo'_{i+\ar{a^1}} (\ul{\sigma},\phi_\sigma\oo a^2_{\inp{\sigma}}) +\ \\
+(a_{\out{\sigma}}\oo_i(\mu_{\out{\sigma}}\oo_2 a^2_{\out{\sigma}})) \oo'_i (\ul{\sigma},\phi_\sigma\oo a^1_{\inp{\sigma}}) +\ \\
-a_{\out{\sigma}}\oo'_i( \ul{\sigma},\phi_\sigma\oo\mu_{\inp{\sigma}}\oo(a^1_{\inp{\sigma}},a^2_{\inp{\sigma}}) )
\end{gather*}
for all $a,a^1,a^2\in\oAss$ and $1\leq i \leq \ar{a}$.
\end{sublemma}

\begin{proof}[of Sublemma \ref{SublemmaFUJ}]
There is a morphism $\psi$ of collections
$$a_{\out{\sigma}}\oo'_i(\ul{\sigma},x_\sigma\oo(a^1_{\inp{\sigma}},\ldots,a^{\ar{x}}_{\inp{\sigma}})) \mapsto a_c\oo'_i(\id_c,x_c\oo(a^1_c,\ldots,a^{\ar{x}}_c))$$
for $x\in X$ (or $x=\phi$) and $a,a^1,a^2,\ldots\in\oAss$.
$\psi$ is obviously an isomorphism of degree $-\dg{\sigma}$.
The differential on the suspension is $(-1)^{\dg{\sigma}}\dd^0$, hence we must verify $$\psi\dd^0 = (-1)^{\dg{\sigma}}\dd^0\psi.$$
This is immediate by the formula \eqref{Defdd0} defining $\dd^0$.
\end{proof}

Using the relations $f\oo\mu_{\inp{f}}-\mu_{\out{f}}\oo(f,f)$ in $\oA$ for $f\in\Sigma^1$, every element $a\in\FrM{\oA}{\Phi_\sigma\op X_\sigma}$ can be written in the form 
\begin{gather} \label{Formulaatop}
a=a_{top}\oo_m(a_{mod}\oo(c_1,\ldots,c_{\ar{a_{mod}}}))
\end{gather}
for some $a_{top}\in\oA$, $m\in\N$, $a_{mod}\in\oM_\sigma$ and $c_1,c_2,\ldots\in\oC$.
We want to make this is expression as close to being unique as possible.
We will require:

\begin{quote}
If $a_{top}$ can be written in the form $$a_{top}=a'\oo_k(\mu_{\out{\sigma}}\oo(a'',\id_{\out{\sigma}})) \quad\mbox{or}\quad a_{top}=a'\oo_k(\mu_{\out{\sigma}}\oo(\id_{\out{\sigma}},a''))$$
for some $a',a''\in\oA$, $1\leq k\leq\ar{a'}$ satisfying $m=k+\ar{a''}$ resp. $m=k$, then $\mu_{\out{\sigma}}\oo(a'',a_{mod})$ resp. $\mu_{\out{\sigma}}\oo(a_{mod},a'')$ can't be written as an element of $\oM_\sigma\oo\oC$.
\end{quote}

It is easily seen that the above requirement can be met for any $a$ so assume it holds in the above expression \eqref{Formulaatop}.
It is easy to see that this determines $a_{top}$ and $a_{mod}$ uniquely up to scalar multiples.
The elements $c_1,c_2,\ldots\in\oC$ are however not unique as the following example shows:

Let $x\in X$ and let $f,g_1,g_2\in\oC$ be such that $fg_1=fg_2$, hence $$a:=x_{\out{f}}\oo(fg_1,\phi_f)=x_{\out{f}}\oo(fg_2,\phi_f).$$
In "the canonical" form \eqref{Formulaatop}, $a_{mod}=x_{\out{f}}\oo(f,\phi_f)$, however $c_1$ is either $g_1$ or $g_2$, in pictures:
$$
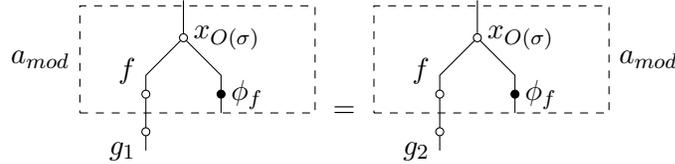
\raisebox{-0.7cm}{
\begin{tikzpicture}[scale=0.5]
\draw (0,1) -- (0,0) node[right]{$x_{\out{\sigma}}$} -- (-1,-1) -- (-1,-1.5) node[above left]{$f$}  -- (-1,-2.5) node[below left]{$g_1$} -- (-1,-3);
\draw (0,0) -- (1,-1) -- (1,-1.5) node[right]{$\phi_f$} -- (1,-2) ;
\draw[dashed] (-2.75,0.85) rectangle (3.5,-2);
\draw (-2.75,-0.575) node[left]{$a_{mod}$} -- cycle;
\filldraw[fill=white] (0,0) circle (3pt);
\filldraw[fill=white] (-1,-1.5) circle (3pt);
\filldraw[fill=white] (-1,-2.5) circle (3pt);
\filldraw (1,-1.5) circle (3pt);
\end{tikzpicture}
}
=
\raisebox{-0.7cm}{
\begin{tikzpicture}[scale=0.5]
\draw (0,1) -- (0,0) node[right]{$x_{\out{\sigma}}$} -- (-1,-1) -- (-1,-1.5) node[above left]{$f$}  -- (-1,-2.5) node[below left]{$g_2$} -- (-1,-3);
\draw (0,0) -- (1,-1) -- (1,-1.5) node[right]{$\phi_f$} -- (1,-2) ;
\draw[dashed] (-2.75,0.85) rectangle (3.5,-2);
\draw (3.5,-0.575) node[right]{$a_{mod}$} -- cycle;
\filldraw[fill=white] (0,0) circle (3pt);
\filldraw[fill=white] (-1,-1.5) circle (3pt);
\filldraw[fill=white] (-1,-2.5) circle (3pt);
\filldraw (1,-1.5) circle (3pt);
\end{tikzpicture}
}$$
\bigskip

So far we have shown that there are, for $m\geq 1$, \coll{\Sigma^0}s $\oA^{top}_{\sigma,m}$ (whose description is implicit in the above discussion) such that there is a \coll{\Sigma^0} isomorphism
\begin{gather} 
\FrM{\oA}{\Phi_\sigma\op X_\sigma}\binom{c}{c_1,\ldots,c_N} \cong \label{HugeIso} \\
\bigoplus_{\substack{m,n\geq 1, \\ d\in\Sigma^0}} \oA^{top}_{\sigma,m}\binom{c}{c_1,\ldots,c_{m-1},d,c_{m+n},\ldots,c_N}\ot\left(\frac{\oM_\sigma\oo\oC}{\mathfrak{L}}\right)\binom{d}{c_m\ldots,c_{m+n-1}}, \nonumber
\end{gather}
where $\mathfrak{L}$ is the sub-\coll{\Sigma^0} of $\oM_\sigma\oo\oC$ describing the non-uniqueness mentioned above.
More precisely, it consists of elements
$$\underbrace{(b\oo'_n(\ul{\sigma},x\oo(b_1,\ldots,b_{\ar{x}})))}_{a_{mod}}\oo(f_1,\ldots,f_{\ar{a_{mod}}})$$
for all $a_{mod}\in\oM_\sigma$, $f_1,f_2,\ldots\in\oC$ satisfying that for some $1\leq i\leq n-1$ or $n+\sum_{j=1}^{\ar{x}}\ar{b_j}\leq i\leq\ar{a_{mod}}$ we have $\ul{\sigma}f_i=0$.
The condition $\ul{\sigma}f_i=0$ means that $f_i=f-g$ (up to a scalar multiple) for some $f,g\in\Sigma^1$ and $\ul{\sigma}f=\ul{\sigma}g$.

At this point we suggest the reader to go through the above discussion in the case $\dg{\sigma}=0$ as many things simplify substantially, e.g. $\mathfrak{L}=0$ and we are essentially done by applying the usual K\"unneth formula to \eqref{HugeIso} and then using Sublemma \ref{SublemmaFUJ} which shows
\begin{gather} \label{NoQutient}
H_*(\oM_\sigma\oo\oC,\dd^0) \cong \frac{(\oAss_{\out{\sigma}}\oo'(\kspan{\ul{\sigma}},\Phi_\sigma\oo\oAss_{\inp{\sigma}}))\oo\oC} {\mathfrak{J}'_\sigma},
\end{gather}
where $\mathfrak{J}'_\sigma$ is an analogue of $\mathfrak{J}_\sigma$ above.
For general $\sigma$, we have to get rid of the quotient:

\begin{sublemma} \label{sublemma2}
$$H_*\left(\frac{\oM_\sigma\oo\oC}{\mathfrak{L}},\dd^0\right) \cong \frac{(\oAss_{\out{\sigma}}\oo'(I,\Phi_\sigma\oo\oAss_{\inp{\sigma}}))\oo\oC} {\mathfrak{J}''_\sigma},$$
\end{sublemma}
where again $\mathfrak{J}''_\sigma$ is the corresponding analogue of $\mathfrak{J}_\sigma$.

\begin{proof}[of Sublemma \ref{sublemma2}]
Denote $$\oM_\sigma\oo\oC\xrightarrow{\proj}\frac{\oM_\sigma\oo\oC}{\mathfrak{L}}=:Q$$ the natural projection.
The differential $\dd^0_Q$ on $Q$ inherited from $\FrM{\oA}{\Phi_\sigma\op X_\sigma}$ is given, for $\alpha\in Q$, by $$\dd^0_Q\alpha = \proj\dd^0\tilde{\alpha},$$
where $\tilde{\alpha}\in\oM_\sigma\oo\oC$ denotes any element such that $\proj\tilde{\alpha}=\alpha$.

By \eqref{NoQutient}, to prove the sublemma it suffices to show
\begin{enumerate}
\item $\proj\Ker\dd^0 = \Ker\dd^0_Q$,
\item $\proj\Im\dd^0 = \Im\dd^0_Q$.
\end{enumerate}

For $1.$, let $\alpha\in Q$, $\dd^0_Q\alpha=0$ and we will show there is $\beta\in\oM_\sigma\oo\oC$ satisfying $\dd^0\beta=0$ and $\proj\beta=\alpha$.
Let $\proj\tilde{\alpha}=\alpha$ and let
$$\tilde{\alpha} = \sum_{i\in I}a^i_{mod}\oo(a^i_1,\ldots,a^i_{\ar{a^i_{mod}}})$$
for some index set $I$, $a^i_{mod}\in\oM_\sigma$ and $a^i_j\in\oC$, $i\in I$ and $1\leq j \leq \ar{a^i_{mod}}$,
such that any two ordered $\ar{a^i_{mod}}$-tuples $(a^i_1,\ldots,a^i_{\ar{a^i_{mod}}})$ are distinct for any two distinct $i$'s.
Let $I_1\subset I$ be the set of $i$'s such that $\dd^0 a^i_{mod}=0$ and let $I_2:=I-I_1$.
By our assumption, $\proj\dd^0\tilde{\alpha} = \dd^0_Q\alpha = 0$, hence
$$\dd^0\tilde{\alpha} = \sum_{i\in I_2}(\dd^0 a^i_{mod})\oo(a^i_1,\ldots,a^i_{\ar{a^i_{mod}}}) \in \mathfrak{L}.$$
Thus for every $i$ there is $j$ such that $\ul{\sigma}\oo a^i_j=0$ (by the definition of $\mathfrak{L}$).
Because $\dd^0$ doesn't change $a^i_j$'s, we get
$\sum_{i\in I_2}a^i_{mod}\oo(a^i_1,\ldots,a^i_{\ar{a^i_{mod}}}) \in \mathfrak{L}$ and we set 
$$\beta := \sum_{i\in I_1}a^i_{mod}\oo(a^i_1,\ldots,a^i_{\ar{a^i_{mod}}}).$$
Then obviously $\dd^0\beta=0$ and $\alpha=\proj\beta$, so we have obtained $\proj\Ker\dd^0 \supset \Ker\dd^0_Q$.
The opposite inclusion is obvious.
Also $2.$ is easy.
\end{proof}

Now apply the sublemma \ref{sublemma2} to \eqref{HugeIso}.
This concludes the proof of Lemma \ref{LemmaA}.
\end{proof}

\begin{proof}[of Lemma \ref{LemmaB}]
This is just a straightforward application of Lemma \ref{LemmaAboutrolDR} to operad $\oA:=\oC$ and its bar-cobar resolution $\oR:=\BC{\oC}$ (e.g. \cite{Vallette}).
To see this, recall that $\BC{\oC}$ is a quasi-free \colop{\Sigma^0} generated by \coll{\Sigma^0} $\kspan{\Sigma^{\geq 1}}$, where the degree of $\sigma\in\Sigma^{\geq 1}$ is $\dg{\sigma}-1$.
The derivation differential is given by
\begin{align*}
\dd\underbrace{(\xleftarrow{f_n}\cdots\xleftarrow{f_1})}_\sigma &:= \sum_{i=1}^{n-1}(-1)^{i+n+1}(\xleftarrow{f_n}\cdots\xleftarrow{f_{i+1}})\oo(\xleftarrow{f_i}\cdots\xleftarrow{f_1})\ + \\
&\phantom{:=}+\sum_{i=1}^{n-1}(-1)^{n-i}\underbrace{(\xleftarrow{f_n}\cdots\xleftarrow{f_{i+1}f_i}\cdots\xleftarrow{f_1})}_{\sigma_i}.
\end{align*}
The projection $\BC{\oC}\xrightarrow{\rR}\BC^1{\oC}\cong\oC$ onto the sub-\coll{\Sigma^0} consisting of single generators is a quism.

Then $\olMDR$ of Lemma \ref{LemmaAboutrolDR} is $$\mathcal{M}\overline{\mathcal{D}\BC{\oC}}=\left( \FrM{\oA}{\bigoplus_{c\in\Sigma^0}\Phi_c\op\bigoplus_{\sigma\in\Sigma^{\geq 1}}\Phi_\sigma},\dd_{\mathcal{M}\overline{\mathcal{D}\BC{\oC}}} \right),$$
where all the symbols have the same meaning as in the previous parts of this paper and it is easily checked that the differential is given by the formula \eqref{FormulaBarCobar}.
This resolves $\mathcal{MDC}$, which is readily seen to be the right-hand side in the statement of Lemma \ref{LemmaB}.
\end{proof}

\end{document}